\documentclass{birkjour}
\usepackage{amsmath, amsthm, amsfonts, mathdots}
\usepackage{amssymb}
\usepackage{amsbsy}
\usepackage{mathrsfs}

\input amssym.def
\input amssym
\usepackage{enumerate}
\usepackage{hyperref}
\hypersetup{
    pdftitle={Memory Estimation of Inverse Operators},    
    pdfauthor={Anatoly Baskakov}{Ilya Krishtal},     
    pdfsubject={Memory Estimation of Inverse Operators},   
    pdfcreator={Ilya Krishtal},   
    pdfproducer={Ilya Krishtal}, 
    pdfkeywords={Banach modules}{Beurling spectrum}{Wiener's lemma}, 
     colorlinks,%
    citecolor=blue,%
    filecolor=black,%
    linkcolor=magenta,%
    urlcolor=black
}

\chardef\bslash=`\\ 

\makeatletter
\def\verbatim{\interlinepenalty\@M \@verbatim
   \leftskip\@totalleftmargin\advance\leftskip2pc
   \frenchspacing\@vobeyspaces \@xverbatim}
\makeatother
\newtheorem{thm}{Theorem}[section]
\newtheorem{cor}[thm]{Corollary}
\newtheorem{lem}[thm]{Lemma}
\newtheorem{prop}[thm]{Proposition}

\newtheorem{ass}[thm]{Assumption}

\newtheorem{con}{Conclusion}[section]
\theoremstyle{definition}
\newtheorem{defn}{Definition}[section]

\theoremstyle{remark}
\newtheorem{rem}{Remark}[section]
\newtheorem{exmp}{Example}[section]

\numberwithin{equation}{section}

\newcommand{\begeq}{\begin {equation}}
\newcommand{\eq}{\end{equation}}
\newcommand{\bs}{\begin {split}}
\newcommand{\es}{\end{split}}
\newcommand{\bp}{\begin {prop}}
\newcommand{\ep}{\end {prop}}
\newcommand{\bt}{\begin {thm}}
\newcommand{\et}{\end {thm}}
\newcommand{\bcon}{\begin {con}}
\newcommand{\econ}{\end {con}}
\newcommand{\bc}{\begin {cor}}
\newcommand{\ec}{\end {cor}}
\newcommand{\bl}{\begin {lem}}
\newcommand{\el}{\end {lem}}
\newcommand{\bpf}{\begin {proof}}
\newcommand{\epf}{\end {proof}}
\newcommand{\bi}{\begin {itemize}}
\newcommand{\ei}{\end {itemize}}
\newcommand{\ben}{\begin {enumerate}}
\newcommand{\een}{\end {enumerate}}
\newcommand{\brem}{\begin {rem}}
\newcommand{\erem}{\end {rem}}

\newcommand{\bd}{\begin {defn}}
\newcommand{\ed}{\end {defn}}

\newcommand{\bex}{\begin {exmp}}
\newcommand{\eex}{\end {exmp}}

\newcommand{\la}{\langle}
\newcommand{\ra}{\rangle}

\newcommand{\TTT}{{\T\kern-.44em \T}}
\newcommand{\tTTT}{\widetilde{\T\kern-.44em \T}}

\newcommand{\wS}{{\widetilde\Sigma}}

\newcommand{\ZZ}{{\mathbb Z}}

\newcommand{\CC}{{\mathbb C}}
\newcommand{\NN}{{\mathbb N}}

\newcommand{\HH}{{\mathcal H}}

\newcommand{\s}{\sigma}
\renewcommand{\l}{\lambda}
\renewcommand{\L}{\Lambda}
\renewcommand{\a}{\alpha}

\newcommand{\M}{\mathcal{M}}

\newcommand{\g}{\gamma}

\newcommand{\T}{\mathcal{T}}

\newcommand{\wP}{\widetilde{P}}
\newcommand{\wcP}{\widetilde{\mathcal P}}

\begin{document}

%
%
%
%
%
%
%
%
%


\title[Spectral Analysis   of Perturbed Operator Matrices]
{
Similarity Techniques in the Spectral \\ Analysis   of Perturbed Operator Matrices} 

\author{Anatoly G. Baskakov}
\address{Department of Applied Mathematics and Mechanics \\ Voronezh State University \\ Voronezh 394693 \\ Russia}
\email{anatbaskakov@yandex.ru}
\thanks{The first author is supported in part by the Ministry of Education and Science of the Russian Federation in the frameworks of the project part of the state work quota (Project No 1.3464.2017/4.6).
The second author is supported in part by NSF grant DMS-1322127. The third author is supported in part by RFBR grant 16-01-00197.}



\author{Ilya A. Krishtal}
\address{Department of Mathematical Sciences\\ Northern Illinois University\\ DeKalb, IL 60115 \\ USA}
\email{ikrishtal@niu.edu}
\author{Natalia B. Uskova}
\address{Department of Higher Mathematics and Mathematical Physical Modeling \\ Voronezh State Technical University \\ Voronezh 394026 \\ Russia}
\email{nat-uskova@mail.ru}


\date{\today }




\begin{abstract}
We develop the method of similar operators to study the spectral properties of unbounded perturbed linear operators that can be represented by matrices of various kinds. The class of operators under consideration includes various differential operators with an involution, such as one-dimensional Dirac operators of a certain type.
\end{abstract}
\maketitle

\section{Introduction}\label{intro}

In this paper, we develop a general approach, called the \emph{method of similar operators}, to study the spectral properties of unbounded perturbed linear operators in a certain large class. The main idea of the method is to construct a similarity transform which would allow one to represent the operator as a block-diagonal matrix. Typically, most of the resulting blocks  have a  structure that is easy to analyze and compute the spectral estimates from. For example, in some cases the blocks end up being of very small dimension \cite{BDS11, BKU19}. This makes the method of similar operators a very effective tool for spectral analysis of various classes of differential and difference operators.  In our recent research \cite{BKU19, BKU18}, as well as in \cite[and references therein]{BDS11, BKR17, BP17}, the method was developed for Dirac operators, Hill operators, differential operators with an involution, etc. In each case, a very similar approach and calculations were used. And yet, there were subtle but crucial differences which necessitated 
essential tweaking of the method in all of those instances. The abstract method of similar operators \cite{B94} was found to be too general to be used without a considerable effort. The purpose of this paper is to develop a version of the method of similar operators that is, on one hand, general enough to be used in cases such as \cite{BDS11, BKR17, BKU19, BKU18, BP17} and, on the other hand, specific enough to make the derivation of the above results clear and transparent.
We fully expect the method in this paper to be easily usable for other classes of differential operators as well. We illustrate this with a new example in Section \ref{newex}.

Let us now present a brief description of our setting. More details will be provided in the following section.

Let $H$ be a separable complex Hilbert space and $B(H)$ be the $C^*$-algebra of all bounded linear operators in $H$. By $A: D(A)\subseteq H\to H$ we denote a closed, densely defined linear operator. We  assume that the operator $A$ is \emph{normal}, that is, the adjoint operator $A^*: D(A^*)\subseteq H\to H$ satisfies $D(A^*) =D(A)$ and $\|A^*x\|=\|Ax\|$ for all $x\in D(A)$ \cite{R91}. We also impose rather stringent conditions on the spectrum $\s(A)$ of the operator $A$ and the resolvent operator $R(\cdot\, ; A): \rho(A) \to B(H)$, $R(\l, A) = (A-\l I)^{-1}$, $\l\in\rho(A)= \CC\setminus\s(A)$. As a part of those conditions, we assume that
\[
\s(A) = \bigcup_{n\in\ZZ} \{\l_n\},
\]
where each $\l_n$ is a semisimple eigenvalue of $A$, that is $AP_n = \l_n P_n$, where
$P_n = P(\{\l_n\}, A)$ is the spectral Riesz projection corresponding to the spectral component
$\s_n = \{\l_n\}$. The full list of conditions on the operator $A$ can be found in Section \ref{condAB}. As a consequence of these conditions, the operators $X: D(X)\subseteq H\to H$ studied in this paper are completely characterized by their operator matrices which will be identified with them and denoted by the same symbol; the entries of the matrices will be the operators $X_{mn} = P_mXP_n$, $m,n \in\ZZ$. In particular, the matrix of the operator $A$ is diagonal and $A_{nn} = \l_nP_n$, $n\in\ZZ$. We shall refer to the operator $A$ as a \emph{free} or \emph{unperturbed} operator.

By \emph{perturbation}, we shall mean a linear operator $B$, which is $A$-bounded (see Definition \ref{abdd}) and satisfies the conditions appearing in Section \ref{condAB}. The operator $A-B$, which is the main object of study in this paper, will be referred to as the \emph{perturbed} operator.

The method of similar operators  yields a similarity transform which gives the similarity between the operator $A-B$ and an operator $A-V$ such that the matrix of $V$ is block diagonal. This allows one, inter alia, to obtain estimates of the spectrum of $A-B$ and prove its generalized spectrality and equiconvergence of spectral decompositions.

Typically, the similarity transform is constructed in two stages; in fact, it is convenient to view
it as a product of two similarity transforms. The first one yields the operator $A - Q$ with $Q$ belonging to  the ideal
of Hilbert-Schmidt operators in $H$ denoted by $\mathfrak{S}_2(H)$\footnote{We shall also make use of the ideal of nuclear operators -- $\mathfrak{S}_1(H)$. Recall that in $\mathfrak{S}_1(H)$ the norm is given by $\|X\|_1=\sum\limits_{n=1}^\infty s_n$, where $(s_n)$ is the
sequence of singular values of the operator $X$. The norm in $\mathfrak{S}_2(H)$  is
$\|X\|_2=(\mathrm{tr}\,XX^*)^\frac{1}{2} =\left(\sum\limits_{n=1}^\infty s_n^2\right)^{1/2}$. We refer to \cite{DS88II, GK65} for the standard properties of these ideals used in this paper.}.
The second one yields the operator $A-V$ with $V$ not only in $\mathfrak{S}_2(H)$ but also having a block diagonal matrix (as mentioned above). Depending on the spectral structure of $A$, we may need to use a smaller space  that contains $Q$ and $V$ and is embedded in $\mathfrak{S}_2(H)$.  For example, a space of this kind was needed in \cite{BKU19, BKU18}, where $A$ was a first order differential operator but not in \cite{BP17}, where $A$ was the second order differential operator. Constructing this smaller space (see Section \ref{wspace}) is the most technical part of this paper.

The remainder of the paper is organized as follows. In Section \ref{prelim}, we introduce the notation used throughout the paper, exhibit conditions commonly assumed for the operators $A$ and $B$, and recall the standard facts about similar operators. Section \ref{method} is devoted to the development of the method of similar operators in the set-up of operator matrices. In Section \ref{main}, we collect our main results on spectral properties of the operator $A-B$ that can be obtained via the method of similar operators. Finally, Section \ref{examples} outlines several examples of such analysis for differential and integro-differential operators. 

\section{Notation and preliminaries}\label{prelim}

In this section, we describe the conditions on the operators $A$ and $B$ that are assumed throughout the paper. We also remind the reader the basic facts about similar operators that are relevant for us and set the notation along the road.

\subsection{Conditions on the operators $A$ and $B$.}\label{condAB}

The list of conditions on the operator $A$ is as follows.

\begin{ass}\label{assa}
The operator $A$ is assumed to have the following properties.
\ben
\item\label{assa1} $A: D(A)\subseteq H\to H$ is a normal closed linear operator. 
\item\label{assa2} The spectrum $\s(A)$ of the operator $A$ satisfies
\[
\s(A) = \bigcup_{n\in\ZZ} \{\l_n\},
\]
where each $\l_n$ is a semisimple eigenvalue of $A$ of finite multiplicity.\footnote{The finite multiplicity condition is not necessary for developing the theory. 
The exposition, however, would become too cumbersome if it is not imposed.}
\item\label{assa3} The operator $A$ is spectral, i.e. the Riesz projections $P_n = P(\{\l_n\}, A)$ form a resolution of the identity.
\item\label{assa4} The eigenvalues $\{\l_n: n\in\ZZ\}$ of $A$ satisfy
\begeq\label{specAcond}
\eta =\sup_{j\in\ZZ}\sum_{n\in\ZZ\setminus\{j\}} |\l_n-\l_j|^{-2}< \infty.
\eq
\item\label{assa5} For any $\varepsilon > 0$ there is $\lambda_\varepsilon\in\rho(A)$,
such that $\|(A-\lambda_\varepsilon I)^{-1}\|<\varepsilon$.
\een
\end{ass} 

 As a consequence of \eqref{specAcond}, we get that the spectrum $\s(A)$ is a separated set, that is 
\begeq\label{speccond}
\delta =\inf\{|\l_m -\l_n|: m\neq n\in\ZZ\} > 0.
\eq
It also follows from \eqref{specAcond} and the normality of $A$ that Assumption \ref{assa}\eqref{assa5} is redundant. It is, however, convenient for our exposition to include it explicitly. 


To set the conditions on the operator $B$, 
we will need the following definition.

\bd\label{abdd}
 Let $A: D(A) \subset H \to H$ be a linear operator. A linear operator $B: D(B) \subset H \to H$ is $A$-bounded if $D(B) \supseteq D(A)$ and $\|B\|_A = \inf\{c > 0: \|Bx\| \le c(\|x\| + \|Ax\|),\ x \in D(A)\} < \infty$.
\ed
The space $\mathfrak L_A(H)$ of all $A$-bounded linear operators with the  domain equal to $D(A)$ is a Banach space with respect to the norm $\|\cdot\|_A$. For densely defined operators $A$, restricting the domain of bounded operators to $D(A)$, allows us to view $B(H)$ as a subspace of $\mathfrak L_A(H)$. 

Recall that, given $\l_0 \in \rho(A)$, where $\rho(A) = \CC\backslash\s(A)$ is the resolvent set of $A$, we have $B \in \mathfrak L_A(H)$ if and
only if $B(\l_0I - A)^{-1} \in B(H)$ and $\|B\|_{\l_0} = \|B(\l_0I - A)^{-1}\|_{B(H)}$ defines an equivalent norm in $\mathfrak L_A(H)$ \cite{EN00}.

In Section \ref{main}, we provide various results for operators $A-B$, where the perturbation $B$ satisfies some combination of the assumptions collected below. Here we formulate the assumptions in terms of the matrix of the operator $B$. For some of the conditions, an equivalent formulation will be given in Section \ref{main}. 
In some cases, weaker versions of the conditions below will be used.

The following assumption will commonly be imposed.

\begin{ass}\label{asB}
The operator $B$ is assumed to have the following properties.
\ben
\item\label{assb1} $B\in \mathfrak L_A(H)$. 
\item\label{assb2} The matrix elements of the operator $B$ satisfy
\begeq\label{cGB}
\sum_{m\in\ZZ}\sum_{n\neq m\in\ZZ} \frac{\|B_{mn}\|^2_2}{|\l_m-\l_n|^2}<\infty,
\eq
and
\begeq\label{cBGB}
\sum_{m,n\in\ZZ}\left\|\sum_{\ell\in\ZZ\setminus\{n\}} \frac{B_{m\ell}B_{\ell n}}{\l_\ell-\l_n}\right\|^2_2<\infty,
\eq
where  $\l_n$, $n\in\ZZ$, are  eigenvalues of $A$.
\item\label{assb3} For any $\varepsilon > 0$ there is $\lambda_\varepsilon\in\rho(A)$,
such that the matrix with elements
\begeq\label{BJA}
Y_{mn} =
\begin{cases}
 \frac{1}{\l_n-\l_\varepsilon}B_{mn}, & m\neq n;\\
0, & m=n; 
\end{cases}
\eq
defines an operator  $Y\in B(H)$ with $\|Y\|<\varepsilon$.
\een
\end{ass}

The following assumption is typically fulfilled in applications.

\begin{ass}\label{asBc}
The matrix elements of the operator $B$ satisfy
\begeq\label{cJB}
\sum_{n\in\ZZ}  {\|B_{nn}\|^2_2}<\infty.
\eq
\end{ass}

\brem\label{compB}
Assumptions \ref{asB} and \ref{asBc} automatically hold if $B\in\mathfrak S_2(H)$. Condition \eqref{cJB} is equivalent to the main diagonal of the matrix of $B$ defining an operator in $\mathfrak S_2(H)$.
\erem

In some examples, we will use an assumption that is weaker than \eqref{cJB}.

\begin{ass}\label{asBb}
The main  diagonal of the matrix of the operator $B$ defines an operator in $B(H)$.
\end{ass}

\brem\label{remB}
If Assumption \ref{asBb} holds, Assumption \ref{asB}\eqref{assb3} is equivalent to
\emph{\ben
\item[3$\,^\prime$.]\label{assb3'} For any $\varepsilon > 0$ there is $\lambda_\varepsilon\in\rho(A)$,
such that $\|B(A-\l_\varepsilon)^{-1}\|<\varepsilon$.
\een}
\erem

\subsection{Operators and their block matrices.}\label{bmat}

We proceed with a more thorough description of the relation between the operators studied in this paper and their matrices. As mentioned in the introduction, the matrices are defined by means of a resolution of the identity. We remind the reader the standard definition, where by $\mathcal I \subseteq \ZZ$ we denote a finite or countable index set.

\bd\label{resid}
A family $\mathcal E = \{E_n: n\in\mathcal I\}$ of idempotents in $B(H)$ is called a \emph{(disjunctive)
resolution of the identity} if the following properties hold:
\ben
\item Each $E_n$ is an orthogonal projection and $E_mE_n = 0$ for $m\neq n\in\mathcal I$.
\item For each $x\in H$, the series $\sum_{n\in\mathcal I} E_nx$ converges unconditionally to $x$.
\een
\ed 

The following property of a resolution of the identity is useful.

\bl\label{propresC}
Let $\mathcal E = \{E_n: n\in\mathcal I\}$ be a resolution of the identity and $\{y_n: n\in\mathcal I\}$ be a sequence of vectors in $H$. Then the series $\sum_{n\in \mathcal I} E_n y_n$ converges unconditionally in $H$ if and only if $\sum_{n\in \mathcal I} \|E_n y_n\|^2 < \infty$.
\el

\bpf
To prove unconditional convergence of the series $\sum_{n\in \mathcal I} E_n y_n$, it suffices to show that for any rearrangement of $\mathcal I$ the sequence of partial sums of 
the series  is Cauchy as long as $\sum_{n\in \mathcal I} \|E_n y_n\|^2 < \infty$. This follows from
\[
\left\|\sum_{n\in\mathcal J} E_ny_n\right\|^2 = 
\left\la \sum_{m\in\mathcal J} E_my_m, \sum_{n\in\mathcal J} E_ny_n\right\ra
=\sum_{n\in\mathcal J} \|E_ny_n\|^2,
\]
where $\mathcal J$ is any (nonempty) finite subset of $\mathcal I$. 

The converse direction follows trivially from Definition \ref{resid}.
\epf

In this paper, we consider only the resolutions of the identity that consist of the spectral projections of the operator $A$. The basic resolution of the identity $\mathcal P = \{P_n: n\in\ZZ\}$ was already mentioned in the introduction. We will, however, also need coarser resolutions that are based on general partitions $\Sigma = \{\s_n: n\in\mathcal I\}$ of the spectrum $\s(A)$. 

\bd\label{cpart}
Let $\Sigma = \{\s_n: n\in\mathcal I\}$ be a partition of the spectrum $\s(A)$ and 
$$P_{\s_n} = \sum_{n: \l_n\in \s_n} P_n= \sum_{\l_n\in \s_n} P(\{\l_n\}, A) = P(\s_n, A).$$ 
We call $\mathcal P_\Sigma = \{P_{\s_n}: n\in \mathcal I\}$ the \emph{$\Sigma$-resolution of the identity}.
\ed

\bex\label{pex}
Besides the basic resolition $\mathcal P$, two kinds of $\Sigma$-resolutions are especially important for us and deserve a special notation. Given $m\in\ZZ$, we let $\mathcal P_m = \{P_m, I-P_m\}$ and
$\mathcal P_{(m)} = \{P_{(m)}\}\cup\{P_n: |n|> m\}$, where
\[
P_{(m)} = \sum_{|n|\le m} P_n.
\]
\eex

In this paper, we make a distinction between operator matrices and matrices of operators.  The definitions are as follows.

\bd
Let $\mathcal P_\Sigma= \{P_{\s_n}: n\in\mathcal I\}$ be a $\Sigma$-resolution of the identity. By a \emph{$\Sigma$-matrix} $X^\Sigma = (X^\Sigma_{mn})_{m,n\in\mathcal I}$ we mean a map $X^\Sigma:\mathcal I \times\mathcal I \to \mathfrak L_A(H)$ such that $X^\Sigma_{mn} = P_{\s_m}X^\Sigma_{mn}P_{\s_n}$. The vector space of all $\Sigma$-matrices will be denoted by $\mathfrak L_A(H, \mathcal P_\Sigma)$.
\ed 

\bd
Let $X\in \mathfrak L_A(H)$ and $\mathcal P_\Sigma= \{P_{\s_n}: n\in\mathcal I\}$ be a $\Sigma$-resolution of the identity. By the \emph{$\Sigma$-matrix of the operator $X$}, we mean the matrix $X^\Sigma$ defined by $X^\Sigma_{mn} = P_{\s_m}XP_{\s_n}$.
\ed

Observe that the $\Sigma$-matrix of an operator is well defined for any $X\in \mathfrak L_A(H)$. Moreover, if $\Sigma$ is such that every ${\s_n}$ is finite, then each $X^\Sigma_{mn} \in B(H)$. By identifying the operators with their matrices, we view $\mathfrak L_A(H)$ as a subspace of $\mathfrak L_A(H, P_\Sigma)$. Conversely, given a matrix $X^\Sigma\in \mathfrak L(H, P_\Sigma)$, we define an operator $X: D(X)\subseteq H\to H$ in the following way. We have $x\in D(X)\subseteq D(A)$ if 
$\sum_{m,n\in\mathcal I}  X^\Sigma_{mn} x$ converges unconditionally to $y\in H$. For $x\in D(X)$,
we let
$Xx= \sum_{m,n\in\mathcal I} X^\Sigma_{mn} x
$.
It may, however, happen that 
an operator defined in this way does not belong to $\mathfrak L_A(H)$. This occurs, for example, if $\mathcal P_\Sigma  = \mathcal P$ and $X^\Sigma_{mn} = n\l_nP_{m}P_{n}$.
Thus, $\mathfrak L_A(H)\subsetneq\mathfrak L_A(H, \mathcal P_\Sigma)$. Observe also that if $\mathcal P_{\Sigma_1}$ is a coarser resolution of the identity than $\mathcal P_{\Sigma_2}$ we may regard $\mathfrak L_A(H, \mathcal P_{\Sigma_1})$ as a subspace of $\mathfrak L_A(H, \mathcal P_{\Sigma_2})$, again, via an appropriate identification of operators and their matrices. 
In general, it may be difficult to determine when matrices in $\mathfrak L_A(H, P_\Sigma)$ come from operators in $\mathfrak L_A(H)$ or $B(H)$. For example, if 
\[
X_{mn} = \begin{cases}
 \frac{1}{m-n}, & m\neq n;\\
0, & m=n; 
\end{cases}\quad\mbox{and}\quad
Y_{mn} =\begin{cases}
 \frac{1}{m-n}, & m > n;\\
0, & m\le n; 
\end{cases}, \ m,n\in\ZZ,
\]
then $X\in B(\ell^2(\ZZ))$ but $Y\notin B(\ell^2(\ZZ))$.
There is, however, a class of operators for which a relatively simple condition exists. An analog of the definition below appeared in \cite{BKU18}.

\bd\label{HScD}
Let $\mathcal E$ be a resolution of the identity.
We say that an operator $X\in B({H})$ belongs to 
the \emph{Hilbert--Schmidt class $\mathfrak{S}_2({H}, \mathcal{E})$
with respect to  $\mathcal{E}$},
if
\begeq\label{Fn}
\sum_{m,n\in\mathcal{I}}\|{E}_m X{E}_n\|^2<\infty.
\eq
\ed

The norm $\|X\|_{2, \mathcal{E}}=\Big(\sum\limits_{m, n\in\mathcal{I}}\|{E}_mX{E}_n\|^2\Big)^{1/2}$ turns $\mathfrak{S}_2({H}, \mathcal{E})$ into a normed linear space.
Moreover, the following three lemmas are immediate.

\bl\label{basklh3}
The space $\mathfrak{S}_2({H}, \mathcal{E})$ of Hilbert--Schmidt operators with respect to a resolution of the identity $\mathcal E$ is a Banach algebra.
\el

\bl\label{basklh5}
An operator $X\in\mathfrak{S}_2({H}, \mathcal{E})$ belongs to $\mathfrak{S}_2({H})$ if and only if
$$
\sum_{m, n\in\mathcal{I}}\|{E}_mX{E}_n\|_2^2<\infty.
$$
\el

\bl\label{lhmo}
Assume that $X^\Sigma\in \mathfrak L_A(H, P_\Sigma)$ is such that $$\sum_{m, n\in\mathcal{I}}\|X^\Sigma_{mn}\|_{B(H)}^2<\infty.$$ Then there is $X\in \mathfrak{S}_2({H}, \mathcal{P}_\Sigma)$ such that $X^\Sigma$ is the $\Sigma$-matrix of $X$.
\el

\brem\label{basklh4}
Observe that  if there is an $N\in\NN$ such that for each $\ell\in\ZZ$ the rank of 
${P}_\ell\in \mathcal P$ is at most $N$, then $\mathfrak{S}_2({H}, \mathcal{P})=\mathfrak{S}_2({H})$. We also note that for any $m\in\ZZ$ we have
$\mathfrak{S}_2({H}, \mathcal{P}_{(m)})=\mathfrak{S}_2({H},\mathcal{P})$ 
and $\mathfrak{S}_2({H}, \mathcal{P}_{m})$ consists of all operators in $B(H)$;
recall that the families $\mathcal P_{(m)}$ and $P_m$ were defined in Example \ref{pex}.
\erem

We conclude this subsection with the following useful result.

\bl\label{lemspecsum}
Assume that an operator $X\in \mathfrak L_A(H)$ has compact resolvent and the $\Sigma$-matrix of $X$ is diagonal. Then
\begeq\label{specsum}
\s(X) = \bigcup_{n\in\mathcal I} \s\left(X^\Sigma_{nn}\vert_{P_{\s_n}H}\right),
\eq
where $X^\Sigma_{nn}\vert_{P_{\s_n}H}$ is the restriction of $X^\Sigma_{nn}$ to the range of the projection $P_{\s_n}$.
\el

\bpf
The ``$\supseteq$'' containment of \eqref{specsum} is trivial. We need to prove the opposite one.
Assume $\l\in\s(X)$. Since the resolvent of $X$ is compact, we have that there exists $x\in D(X)$ such that $Xx = \l x$. Since the matrix $X^\Sigma$ is diagonal, we have $XP_{\s_n}x = \l P_{s_n}x$ for all $n\in\mathcal I$. Since $x\neq 0$ and $\mathcal P_\Sigma$ is a resolution of the identity, we must have $\mathcal P_{\Sigma_n}x\neq 0$ for some $n\in\mathcal I$. For that $n$ we clearly have that 
$P_{\Sigma_n}x$ is an eigenvector of $X^\Sigma_{nn}\vert_{P_{\s_n}H}$ corresponding to the eigenvalue $\l$, and the result follows.
\epf

\brem
We remark that \eqref{specsum} is not true in general, not even in the case when its right-hand-side is a closed set. For example, if $P_{\s_n}H\simeq\CC^n$ and $X^\Sigma_{nn}\vert_{P_{\s_n}H}$ is the $n\times n$ Jordan cell:
\[
X^\Sigma_{nn}\vert_{P_{\s_n}H} = \left( \begin{array}{cc}
0 & 0 \\
I_{n-1}& 0
\end{array} 
\right) = \left( \begin{array}{cccccc}
0 & 0 &\cdots & 0 &0&0\\
1& 0&\cdots& 0&0&0\\
0& 1&\cdots& 0&0&0\\
\vdots & \vdots & \ddots&\vdots&\vdots&\vdots\\
0&0&\cdots&1&0&0 \\
0&0&\cdots&0&1&0
\end{array} 
\right),
\]
then the union in the right-hand-side is $\{0\}$ whereas the spectral radius $r(X)$ of the operator $X$ is $1$ by the Gelfand formula.
\erem


%
%

\subsection{Basic facts about similar operators.}

\bd\label{defsim}
Two linear operators $A_m: D(A_m)\subset\mathcal{H}\to\mathcal{H}$, $m=1, 2$, are called 
\emph{similar}, if there exists a continuously invertible operator
 $U\in B(\mathcal{H})$ such that
$$
A_1Ux=UA_2x, \quad x\in D(A_2), \quad UD(A_2)=D(A_1).
$$
The operator $U$ is called the \emph{similarity transform} of $A_1$ into $A_2$.
\ed

Directly from Definition \ref{defsim}, we have the following result about the spectral properties of similar operators.

\bl\label{basklh1}
Let $A_m: D(A_m) \subset \HH\to\HH$, $m = 1,2$, be two similar operators with the operator $U$ being the similarity transform of $A_1$ into $A_2$. Then the following properties hold.
\ben
\item[(1)]	We have $\s(A_1) = \s(A_2)$, $\s_p(A_1) = \s_p(A_2)$, and $\s_c(A_1) = \s_c(A_2)$, where $\s_p$ denotes the point spectrum and $\s_c$ denotes the continuous spectrum;
\item[(2)] If $\l$ is an eigenvalue of the operator $A_2$ and $x$ is a corresponding eigenvector, then $y = Ux$ is  an eigenvector of the operator $A_1$ corresponding to the same eigenvalue $\l$.
\item[(3)]	If $P_\s$ is the spectral projection of $A_2$ that corresponds to the spectral component $\s\subseteq\s(A_2) = \s(A_1)$, then $ P^\prime_\s = UP_\s U^{-1}$ is the spectral projection of $A_1$ that corresponds to the same spectral component $\s$. 
\een
\el

\section{The method of similar operators}\label{method}

The method of similar operators has its origins in various similarity and perturbation techniques. Among them, there are the classical perturbation methods of celestial mechanics, Ljapunov's kinematic
similarity method \cite{GKK96, Lj56, N15},  Friedrichs' method of similar operators that is used in quantum mechanics \cite{F65}, and Turner's method of similar operators  \cite{T65, U04}. 

The method of similar operators has been extensively developed and used for various classes of unbounded linear operators, see e.g. \cite{B83, B85, B86, B94, B99, B15, BDS11, BK88, BK13, BKR17, BKU19, BKU18}. In this paper, we create a version of the method that can be used to obtain results such as in  \cite{BDS11, BKR17, BKU19, BKU18, BP17}. In Subsection \ref{basmeth}, we exhibit the basic ideas and theorems of the method. In the following Subsection \ref{genJG}, we 
provide a basic construction for 
operator matrices.
In Subsection \ref{coarseJG}, we describe a more elaborate construction of the transforms that is based on a more general resolution of the identity. Finally, in Subsection \ref{wspace}, we construct a weighted space of perturbations that is tailored for a specific perturbation at hand. 

\subsection{Abstract scheme and results.}\label{basmeth}

As we mentioned in the introduction, 
the main idea of the method of similar operators is to construct a similarity transform for the operator $A - B:  D(A) \subset \HH \to \HH$, where   $A$  and  $B$ satisfy the assumptions in Subsection \ref{condAB}. The goal of the method is to obtain an operator $V$ such that   $A-B$ is similar to $A-V$ and  the spectral properties of $A-V$ are in some sense close to those of $A$. In particular, certain spectral subspaces of $A$ are mapped by the similarity transform onto certain subspaces that are invariant for $A-V$.

The method of similar operators uses the \emph{commutator transform} 
$\mathrm{ad}_A:
D(\mathrm{ad}_A)\subset B(\mathcal{H})\to B(\mathcal{H})$ defined by
\begeq\label{ct}
\mathrm{ad}_AX = AX-XA, \quad X\in D(\mathrm{ad}_A).
\eq
The domain $D(\mathrm{ad}_A)$ in \eqref{ct} consists of all  $X\in B(\mathcal{H})$ such that
the following two properties hold:
\ben
\item $XD(A)\subseteq D(A)$;
\item The operator $\mathrm{ad}_AX: D(A)\to\mathcal{H}$ admits a unique extension to a bounded operator $Y\in B(\mathcal{H})$; we then let
$\mathrm{ad}_AX=Y$.
\een

The key notion of the method of similar operators is that of an admissible triplet. Once such a triplet is constructed, achieving the goal of the method becomes a  routine task.

\bd[\cite{BDS11, BP17}]\label{baskdef8}
Let $A$ be a closed, densely defined linear operator, $\mathcal{M}$ be a linear subspace of $\mathfrak{L}_A(\mathcal{H})$,
$J: \mathcal{M}\to\mathcal{M}$, and $\Gamma: \mathcal{M}\to B(\mathcal{H})$.
The collection $(\mathcal{M}, J, \Gamma)$ is  an \emph{admissible triplet} for the operator $A$, and the space
$\mathcal{M}$ is the \emph{space of admissible perturbations}, if the following six properties hold.

\ben
\item\label{adprop1}
 $\mathcal{M}$ is a Banach space that is continuously embedded in $\mathfrak{L}_A(\mathcal{H})$, i.e., $\mathcal{M}$ has a norm $\|\cdot\|_{\mathcal M}$
such that there is a constant $C>0$ that yields $\|X\|_A\le  C\|X\|_{\mathcal M}$ for any
$X\in\mathcal{M}$.

\item\label{adprop2}
 $J$ and $\Gamma$ are bounded linear operators; moreover, $J$ is an idempotent.

\item\label{keyprop} $(\Gamma X)D(A)\subset D(A)$ and
$$
(\mathrm{ad}_A\,\Gamma X)x = (X-JX)x, \quad x\in D(A), \quad   X\in\mathcal{M};
$$
moreover $Y = \Gamma X\in B(\mathcal{H})$ is the unique solution of the equation
\begin{equation}\label{bask11'}
\mathrm{ad}_A\,Y = AY-YA = X-JX,
\end{equation}
that satisfies $JY=0$.

\item\label{adprop4} $X\Gamma Y$, $(\Gamma X)Y\in\mathcal{M}$ for all $X, Y\in\mathcal{M}$, and there is a constant $\gamma>0$ such that
$$
\|\Gamma\|\le\gamma, \quad \max\{\|X\Gamma Y\|_{\mathcal M}, \|(\Gamma X)Y\|_{\mathcal M}\}\le \gamma\|X\|_{\mathcal M}\|Y\|_{\mathcal M}.
$$

\item\label{adprop5} $J((\Gamma X)JY)=0$ for all $X, Y\in\mathcal{M}$.

\item\label{adprop6} For every $X\in\mathcal{M}$ and $\varepsilon>0$ there exists a number  $\lambda_\varepsilon\in\rho(A)$,
such that $\|X(A-\lambda_\varepsilon I)^{-1}\|<\varepsilon$.
\een
\ed

To get an intuitive idea about  the above definition, one should think of the operators involved in terms of their matrices. Recall that the operator $A$ is  represented by a diagonal matrix. The operator $B$ typically has a matrix with some kind of off-diagonal decay. The transform $J$ should be thought of as a projection that picks  the main (block) diagonal of an infinite matrix, whereas the transform $\Gamma$ annihilates the main (block) diagonal and weighs the remaining diagonals in accordance with equation \eqref{bask11'} thereby introducing or enhancing the off-diagonal decay. 
In some vague sense, the method of similar operators includes constructing a sequence of transforms that yield stronger and stronger off-diagonal decay resulting in a block diagonal matrix in the end. 

To formulate the main theorem of the method of similar operators for an operator $A-B$, we use the function
$\Phi:\M\to\M$ given by
\begeq\label{bask13}
\Phi(X) = B\Gamma X-(\Gamma X)(JB)-(\Gamma X)J(B\Gamma X)+B.
\eq 
\bt[\cite{BDS11, BP17}]\label{baskth6}
Assume that $(\mathcal{M}, J, \Gamma)$ is an admissible triplet for an operator $A: D(A)\subset\mathcal{H}\to\mathcal{H}$ 
and $B\in\mathcal{M}$. Assume also that
\begeq\label{bask12}
4\gamma\|J\|\|B\|_{\mathcal M}<1,
\eq
where $\gamma$ comes from  Property \ref{adprop4} of Definition \ref{baskdef8}. Then 
the map $\Phi:\mathcal{M}\to\mathcal{M}$ given by \eqref{bask13} is a contraction and has a unique fixed point 
$X_*$ in the ball 
\begeq\label{ball}
\mathscr B = \{X\in\mathcal{M}: \|X-B\|_{\mathcal M}\le  3\|B\|_{\mathcal M}\}, 
\eq
which can be found as a limit of simple iterations: $X_0=0$, $X_1=\Phi(X_0) = B$, etc.
Moreover,  the 
operator $A-B$ is similar to the operator $A-JX_*$ 
and the similarity transform of 
$A-B$ into $A-JX_*$ is given by $I+\Gamma X_*\in
B(\mathcal{H})$. 
\et

We include 
the proof for completeness.

 \bpf
 We use the Banach fixed-point theorem to prove existence of $X_*$. To apply the theorem, it suffices to show that $\Phi(\mathscr B) \subseteq \mathscr B$ and $\|\Phi(X) -\Phi(Y)\|_\M \le q\|X-Y\|_\M$ for all  $X, Y\in\mathscr B$ and some $q\in(0,1)$.  The first condition follows from \eqref{bask12} and
 \[
 \bs
\|\Phi(X)  - B\|_\M &\le \|B\Gamma X - (\Gamma X)(JB) - (\Gamma X)J(B\Gamma X) \|_\M \\ &
\le 4\g\|B\|_\M^2+ 4\g\|J\|\|B\|_\M^2+16\g^2\|J\|\|B\|_\M^3 \\ &
<  \left(\|J\|^{-1}+1+\|J\|^{-1}\right)\|B\|_\M \le 3\|B\|_\M,
 \end{split}
 \]
 where we used the fact that $ \|X\|_\M \le (\|X-B\|_\M+\|B\|_\M)\le 4\|B\|_\M $ for all $X\in\mathscr B$.

 To obtain the second condition, we let $q = 4\g\|J\|\|B\|_\M <1$ and estimate
 \[
 \bs
 \|\Phi(X)  - \Phi(Y) \|_\M  &\le 
 \|B\Gamma(X-Y)\|_\M+\|(\Gamma(X-Y))JB\|_\M \\&
 +\|(\Gamma X)J(B\Gamma X) - (\Gamma Y)J(B\Gamma Y)\|_\M \\&
 \le \g\|B\|_\M\|X-Y\|_\M + \g\|J\|\|B\|_\M\|X-Y\|_\M \\&
 +\|(\Gamma X)J(B\Gamma (X-Y))\|_\M + \|(\Gamma (X-Y))J(B\Gamma Y)\|_\M \\ &
 \le \g\|B\|_\M\|X-Y\|_\M + \g\|J\|\|B\|_\M\|X-Y\|_\M \\&
 +\g^2\|J\|\|B\|_\M(\|X\|_\M+\|Y\|_\M)\|X-Y\|_\M \\ &
\le \left(\frac q4\|J\|^{-1} +\frac q4 + \frac {q^2}2\|J\|^{-1}\right)\|X-Y\|_\M
\le q\|X-Y\|_\M.
 \end{split}
 \]
 Now the Banach fixed-point theorem implies that the sequence $\{X_n\}$ in the formulation of the theorem converges to an operator $X_*\in\mathscr B\subset\M$ such that $\Phi(X_*) = X_*$. Next, we need to verify that 
 \begeq\label{mainsim}
 (A-B)(I+\Gamma X_*) = (I+\Gamma X_*)(A-JX_*).
 \eq
From  Property \ref{adprop5} of Definition \ref{baskdef8} and \eqref{bask13}, we have
\begeq\label{JXs}
 JX_* = J(B\Gamma X_*)+JB.
 \eq
 Using \eqref{bask11'}, \eqref{bask13} and \eqref{JXs}, we compute
 \[
 \bs
  &(A-B)(I+\Gamma X_*) = A - B +A\Gamma X_* -B\Gamma X_* \\&
  = A-B+(\Gamma X_*)A +X_*-JX_*-B\Gamma X_* \\&
  =A+(\Gamma X_*)A -(\Gamma X_*)JB - (\Gamma X_*)J(B\Gamma X_*)-JX_* \\&
  =A +(\Gamma X_*)A -(\Gamma X_*)JX_*- JX_* = (I+\Gamma X_*) (A-JX_*),
 \end{split}
 \]
 proving \eqref{mainsim}. Next, we need to show that $I+\Gamma X_*$ is continuously invertible. Using $X_*\in\mathscr B$, Property \ref{adprop4}, and \eqref{bask12}, we get
 \begeq\label{GXs}
 \|\Gamma X_*\| \le \g \|X_*\|_\M \le \g(\|X_*-B\|_\M+\|B\|_\M)\le 4\g\|B\|_\M < 1,
 \eq
 yielding 
  \begeq\label{NGX} (I+\Gamma X_*)^{-1} = \sum_{n=0}^\infty (-1)^n (\Gamma X_*)^n.\eq 
  To complete the proof of the theorem, it suffices to show that  $(I+\Gamma X_*)D(A) = D(A)$.
 From  Property \ref{keyprop} of Definition \ref{baskdef8}, we immediately get  $(I+\Gamma X_*)D(A) \subseteq D(A)$. It remains to show that
 $(I+\Gamma X_*)^{-1}D(A) \subseteq D(A)$. Using Property \ref{keyprop} once again, we get
  \begeq\label{bgg}
 \bs
 \Gamma X_*(A &-\l I)^{-1}  = (A-\l I)^{-1}(A-\l I)\Gamma X_* (A-\l I)^{-1} \\
 & = (A-\l I)^{-1}( X_* - JX_* +\Gamma X_* A- \l \Gamma X_*)(A-\l I)^{-1} \\
 & = (A-\l I)^{-1}((X_*-JX_*)(A-\l I)^{-1}+\Gamma X_*)
 \end{split}
 \eq
 for any $\l\in\rho(A)$. Using Property \ref{adprop6} and \eqref{GXs}, we choose $\l\in\rho(A)$ such that
 \[
 \|(X_*-JX_*)(A-\l I)^{-1}+\Gamma X_*\|<1.
 \]
 Finally, using \eqref{NGX} and \eqref{bgg}, we get
 \[
 (I+\Gamma X_*)^{-1} (A-\l I)^{-1} = (A-\l I)^{-1}(I+(X_*-JX_*)(A-\l I)^{-1}+\Gamma X_*)^{-1}
 \]
 and the theorem is proved.
 \epf

\brem\label{rem3}
If the perturbation $B$ in Theorem \ref{baskth6} is such that $JB =0$, then the condition \eqref{bask12} may be weakened to
$3\gamma\|J\|\|B\|_{\mathcal M}<1$. The proof follows the same lines as above.
\erem

There are various appropriate choices for the space $\mathcal M$ in the above theorem. The correct choice usually depends on the perturbation $B$. Sometimes, one can use $B(H)$ or 
$\mathfrak S_2(H)$ in place of $\mathcal M$. In other cases, one is forced to use a narrower space, which we exhibit in Subsection \ref{wspace}.

Often, it can be difficult to pick the correct space $\mathcal{M}$ in advance. Under our assumptions on $A$ and $B$, however, it may be possible to make a preliminary choice of $\mathcal{M}$ first, and then find an operator $A-Q$ that is similar to $A-B$ and such that $Q\in \mathcal{M}$. To state the result, we need the following assumption.

\begin{ass}[\cite{BDS11}]\label{baskpred1}
Assume that $(\mathcal{M}, J, \Gamma)$ is  an {admissible triplet} for an operator $A$
such that the transforms $J$ and $\Gamma$ are restrictions of linear operators from 
 $\mathfrak{L}_A({H},\mathcal P)$ to  $\mathfrak{L}_A({H}, \mathcal P)$ denoted by the same symbols. Assume also that the operator $B\in \mathfrak{L}_A({H})$ has the following five properties.
 
\ben
\item\label{propas1} $\Gamma B\in B(\mathcal{H})$ and $\|\Gamma B\|<1$;

\item\label{propas2} $(\Gamma B)D(A)\subseteq D(A)$;

\item\label{mprop} $B\Gamma B$, $(\Gamma B)JB\in\mathcal{M}$;

\item\label{propas4} $A(\Gamma B)x-(\Gamma B)Ax=Bx-(JB)x$, $x\in D(A)$;

\item\label{propas5} For any  $\varepsilon>0$ there is $\lambda_\varepsilon\in\rho(A)$ such that
$\|(B-JB)(A-\lambda_{\varepsilon}I)^{-1}\|<\varepsilon$.
\een
\end{ass}
\bt[\cite{BDS11}]\label{baskth8}
If Assumption \ref{baskpred1} holds then the operator $A-B$ is similar to   $A-JB-B_0$, where
$B_0=(I+\Gamma B)^{-1}(B\Gamma B-(\Gamma B)JB)$. 
The similarity transform is given by $I+\Gamma B$ so that
\begeq\label{abs1}
(A-B)(I+\Gamma B)=(I+\Gamma B)(A-JB-B_0).
\eq
\et
 
 Since Assumption \ref{baskpred1} is slightly different from its analog in \cite{BDS11}, we include the proof.
 
 \bpf
 
 Condition \ref{propas1} of Assumption \ref{baskpred1}  ensures that the operator $(I+\Gamma B)^{-1}$ is well defined. The operator $B_0$ is then well defined because of Condition \ref{mprop}. 
  Conditions \ref{propas2} and \ref{propas4} yield \eqref{abs1} via a direct computation.
 Since $D(A-B) = D(A-JB-B_0) = D(A)$, it remains to show that $(I+\Gamma B)^{-1}D(A) = D(A)$. The argument is very similar to that  in the proof of Theorem \ref{baskth6}. We use Conditions  \ref{propas2} and \ref{propas4} to obtain
 \[
 \bs
 \Gamma B(A-\l I)^{-1} & = (A-\l I)^{-1}(A-\l I)\Gamma B (A-\l I)^{-1} \\
 & = (A-\l I)^{-1}( B - JB +\Gamma B A- \l \Gamma B)(A-\l I)^{-1} \\
 & = (A-\l I)^{-1}((B-JB)(A-\l I)^{-1}+\Gamma B)
 \end{split}
 \]
 for any $\l\in\rho(A)$. It follows that
 \begeq\label{gneu}
( \Gamma B)^n(A-\l I)^{-1} = (A-\l I)^{-1}((B-JB)(A-\l I)^{-1}+\Gamma B)^n
 \eq
 for any $n\in \NN$. Conditions \ref{propas1} and \ref{propas5} imply that we can choose $\l\in\rho(A)$ such that $\|(B-JB)(A-\l I)^{-1}+\Gamma B\|<1$. Using Neumann series decomposition and \eqref{gneu}, we get
 \[
 (I+ \Gamma B)^{-1}(A-\l I)^{-1} = (A-\l I)^{-1}(I+(B-JB)(A-\l I)^{-1}+\Gamma B)^{-1}
 \]
 yielding $(I+\Gamma B)^{-1}D(A) \subseteq D(A)$. The opposite containment follows from Condition \ref{propas2} and the theorem is proved.
 \epf
 
 In many cases, we have $Q = JB+B_0\in \mathcal M$ in  Theorem \ref{baskth8}. In other interesting cases, we only have  $B_0\in \mathcal M$; then, we will need to apply Theorem \ref{baskth6} for $A-JB$ in place of $A$. 
 
 \brem
 Notice that the assumptions of Subsection \ref{condAB} were not used in the proofs of Theorems \ref{baskth6} and \ref{baskth8}. We only need them to construct specific transforms $J$ and $\Gamma$ in the following subsection.
 \erem
 
\subsection{Standard transforms $J$ and $\Gamma$.}\label{genJG}

As mentioned above, the transform $J$  is supposed to pick out the main diagonal of the operator matrix. Hence, 
we define $J: \mathfrak{L}_A({H},\mathcal P)\to \mathfrak{L}_A({H},\mathcal P)$ via
\begeq\label{tJ2}
(JX)_{mn} = \delta_{m-n} X_{mn},\ m,n\in\ZZ,\  X\in \mathfrak{L}_A({H},\mathcal P),
\eq
where $\delta_{k}$ is the usual Kronecker delta. 

Observe that in this setting the matrix of the commutator  $\mathrm{ad}_AX$ in \eqref{ct} satisfies
\[
(\mathrm{ad}_AX)_{mn} =   (\l_m - \l_n)X_{mn}, \ X\in D(\mathrm{ad}_A).
\]
Therefore, in view of Property \ref{keyprop} of Definition \ref{baskdef8}, it is natural to define the transform $\Gamma$  via
\begeq\label{tG2}
(\Gamma X)_{mn} =
\begin{cases}
 \frac{1}{\l_m-\l_n}X_{mn}, & m\neq n;\\
0, & m=n; 
\end{cases}
\quad X\in \mathfrak{L}_A({H},\mathcal P). 
\eq

The following lemma collects the basic 
properties of the transforms $J$ and $\Gamma$ defined by \eqref{tJ2} and \eqref{tG2}. 
As usually, we identify operators in $\mathfrak{L}_A({H})$ with their matrices.

\bl\label{bask100}
The  following properties hold.
\ben
\item\label{BHext} \label{adprop11}
The transforms $J$ and $\Gamma$  defined by \eqref{tJ2} and \eqref{tG2} respectively restrict to operators in $B(\mathfrak L_A(H))$, $B(B(H))$ and $B(\mathfrak S_2(H,\mathcal P_\Sigma))$, for any partition $\Sigma$ of the spectrum $\s(A)$.  The restrictions will be denoted by the same symbols.
\item\label{adprop21} The transform $J$ is an idempotent and $$\|J\|_{B(B(H))}
= \|J\|_{B(\mathfrak S_2(H,\mathcal P_\Sigma))} = 1.$$
\item\label{adprop31} The transform $\Gamma$   satisfies 
$$\|\Gamma\|_{B(B(H))}\le\sqrt \eta, \|\Gamma\|_{B(\mathfrak S_2(H, \mathcal P_\Sigma))}\le\sqrt \eta, \mbox{ and }  
\|\Gamma\|_{B(\mathfrak S_2(H, \mathcal P))} \le\frac1\delta,$$
where $\eta$ is defined by \eqref{specAcond} and $\delta$ -- by \eqref{speccond}.
\item\label{adprop41}
For  $X\in B(H)$, we have $\Gamma X \in D({\rm{ad}}_{A})$ and
\begeq\label{gax}
{\rm{ad}}_{A}(\Gamma X) =  {A}(\Gamma X) - (\Gamma X){A} = X - J X.
\eq
\een
\el

\bpf
The statements for $\mathfrak S_2(H,\mathcal P)$ follow immediately since we have
\begeq\label{tJ}
JX = \sum_{n\in \ZZ} P_n XP_n, \ X\in \mathfrak S_2(H,\mathcal P), 
\eq
and
\begeq\label{tG}
\Gamma X =\sum_{\substack{m,n \in\ZZ \\ {m\neq n}}}\frac{P_mXP_n}{\l_m-\l_n}, \ X\in \mathfrak S_2(H,\mathcal P), %
\eq
where the series converge unconditionally in $\mathfrak S_2(H,\mathcal P)$. 

The properties $J\in B(B(H))$ and $\|J\|_{B(B(H))} = 1$ follow from Lemma \ref{propresC} and the inequalities
\[
\|JXx\|^2 = \sum_{n\in\mathcal I} \|P_nXP_n x\|^2\le\sum_{n\in\mathcal I}\|X\|^2\|P_n x\|^2=\|X\|^2\|x\|^2, 
\]
which hold for any $x\in H$ and $X\in B(H)$.

Next, we prove that $\Gamma(B(H))\subseteq B(H)$ and $\|\Gamma\|_{B(B(H))}\le\sqrt \eta$.
For $x\in H$, let 
\[
y_{mn} =
\begin{cases}
 \frac{x}{\l_m-\l_n}, & m\neq n;\\
0, & m=n. 
\end{cases}
\] 
Lemma \ref{propresC} and \eqref{speccond} yield $y_m = \sum_{n\in\mathcal I} P_ny_{mn} \in H$ via
\[
\sum_{n\in\mathcal I}\|P_ny_{mn}\|^2 = \sum_{n\in\mathcal I\setminus\{m\}}\left\|\frac{P_ nx}{\l_m-\l_n}\right\|^2 \le \frac{\|x\|^2}{\delta^2} <\infty.
\]
Given $X\in B(H)$, we use Lemma \ref{propresC}  for the sequence $\{Xy_m\}$ and obtain $\|\Gamma X\|_{B(H)} \le \sqrt \eta\|X\|_{B(H)}$ due to
\[\bs
 \left\|(\Gamma X)x\right\|^2 &= \left\|\sum_{m\in\mathcal I} P_mXy_m\right\|^2
 = \sum_{m\in\mathcal I} \left\|P_mXy_m\right\|^2 
\le \|X\|^2  \sum_{m\in\mathcal I} \|y_m\|^2
\\& \le \|X\|^2 \sum_{n\in\mathcal I} \left(\|{P_ nx}\|^2\sum_{m\in\mathcal I\setminus\{n\}} |{\l_m-\l_n}|^{-2} \right)\le \eta\|X\|^2{\|x\|^2}  <\infty.
\end{split}\] 

The statements for $\mathfrak S_2(H, \mathcal P_\Sigma)$ follow from the definition of $\mathfrak S_2(H, \mathcal P_\Sigma)$ and the  above results for $B(H)$.  

Property \ref{adprop41} is now easily verified by direct computation. In particular, one checks that for a given $X\in B(H)$ the bounded operators $(\Gamma X)(A-\l I)^{-1} - (A-\l I)^{-1}\Gamma X$ and $(A-\l I)^{-1}(X-JX)(A-\l I)^{-1}$ have the same matrices. It follows that
$(\Gamma X) D(A) \subseteq D(A)$ and \eqref{gax} is again verified by comparing the matrices.

Finally, the statements for $\mathfrak L_A(H)$ follow from 
\begeq\label{Jext}
JX = J(X(A-\l I)^{-1}) (A-\l I), \quad X\in  \mathfrak L_{A}(H),
\eq
and
\begeq\label{Gext}
\Gamma X = \Gamma (X({A}-\l I)^{-1}) ({A}-\l I), \quad X\in  \mathfrak L_{A}(H),
\eq
which hold for any $\l\in \rho(A)$. 
 Moreover, if we use the $\l$-norm in $\mathfrak L_{A}(H)$, we get
$
\|J\|_{B(\mathfrak L_A(H))} = 1$   and   $\|\Gamma\|_{B(\mathfrak L_A(H))} \le\sqrt \eta.
$
\epf

\brem
In many examples \cite[and references therein]{BDS11, BKU18, BKU19}, the transforms $J$ and $\Gamma$ admit integral representations that yield a better estimate for $\|\Gamma\|_{B(B(H))}$ compared to the one in Lemma \ref{bask100}\eqref{adprop31}.
\erem

An analog of the following lemma for the case $\mathfrak S_2(H, \mathcal P)=\mathfrak S_2(H)$ appears in \cite{BKU19} . 
In Lemma \ref{JG2} below we will prove a more general result.

\bl[{\cite[Lemma 3.5]{BKU19}}]\label{JG1}
Assume that the operator $A$ satisfies Assumption \ref{assa}, $\M = \mathfrak S_2(H, \mathcal P)$, and the transforms $J$ and $\Gamma$ are defined by \eqref{tJ} and \eqref{tG}. 
Then $(\M, J, \Gamma)$ is an admissible triplet for   $A$.
\el

%

\subsection{Coarser transforms $J_\Sigma$ and $\Gamma_\Sigma$.}\label{coarseJG}

Next, we use the above construction to define transforms $J_\Sigma$ and $\Gamma_\Sigma$ for an arbitrary resolution of the identity $\mathcal P_\Sigma$. Naturally, we have
\begeq\label{tJ2s}
(J_\Sigma X)_{mn} = \delta_{m-n} P_{\s_m}X_{mn}P_{\s_n},\ m,n\in\mathcal I,\  X\in \mathfrak{L}_A({H},\mathcal P),
\eq
and
\begeq\label{tJs}
J_\Sigma X = \sum_{n\in \ZZ} P_{\s_n} XP_{\s_n}, \ X\in \mathfrak S_2(H,\mathcal P_\Sigma). 
\eq

To 
enforce Property \ref{keyprop} of Definition \ref{baskdef8}, we  let
\begeq\label{Gs}
\Gamma_\Sigma X = \Gamma (X- J_\Sigma X), \quad X\in\mathfrak{L}_A({H},\mathcal P). %
\eq

\bl\label{JG2}
Assume that the operator $A$ satisfies Assumption \ref{assa}, $\M = \mathfrak S_2(H, \mathcal P)$, $\M_\Sigma =  \mathfrak S_2(H, \mathcal P_\Sigma)$, and the transforms $J_\Sigma$ and $\Gamma_\Sigma$ are defined by \eqref{tJs} and \eqref{Gs}.
Then $(B(H), J_\Sigma, \Gamma_\Sigma)$, $(\M, J_\Sigma, \Gamma_\Sigma)$, and $(\M_\Sigma, J_\Sigma, \Gamma_\Sigma)$ are admissible triplets for $A$. 
\el


\bpf
Property \ref{adprop1}  of Definition \ref{baskdef8}  is immediate from the definitions of $\M$ and $\M_\Sigma$. Property \ref{adprop2} follows from the definitions of $J_\Sigma$ and $\Gamma_\Sigma$, Lemma \ref{propresC}, and Lemma \ref{bask100}(\ref{adprop11}--\ref{adprop31}). Property \ref{keyprop} follows from \eqref{Gs} and Lemma \ref{bask100}(\ref{adprop41}). Indeed,  since 
$(\Gamma X)D(A) \subseteq D(A)$ for all $X\in B(H)$ due to Lemma \ref{bask100}(\ref{adprop41}), we have $(\Gamma_\Sigma X)D(A) \subseteq D(A)$ from \eqref{Gs}. From Lemma \ref{bask100}(\ref{adprop41}) and \eqref{Gs}, we also have
\[
A\Gamma_\Sigma X - (\Gamma_\Sigma X)A = (X-J_\Sigma X) - J(X - J_\Sigma X) = X - J_\Sigma X, \ X\in B(H).
\]
The uniqueness condition of Property \ref{keyprop} follows since the only matrices that commute with $A$ are diagonal.
Property \ref{adprop4} follows from \ref{bask100}(\ref{adprop31}), \eqref{Gs}, and Lemma \ref{basklh3}.
Property \ref{adprop5} is easily verified by direct computation on matrices. It can also be obtained as an  application of \cite[Corollary 7.8]{BK05}. Finally, Property  \ref{adprop6} follows from Assumption \ref{assa}(\ref{assa5}).
\epf


\subsection{The weighted Hilbert-Schmidt spaces.}\label{wspace}

In this section, we let $\Sigma = \{s_n: n\in\ZZ\}$ be some fixed partition of $\s(A)$.
To circumvent condition \eqref{bask12} in Theorem \ref{baskth6} and be able to handle larger perturbations, we sometimes need a weighted version of the space $\M_\Sigma = \mathfrak S_2(H, \mathcal P_\Sigma)$.  In this subsection, we define such a space  $\mathcal M_\Sigma^B \subseteq \M_\Sigma$ for any given $B\in \M_\Sigma$. To simplify the notation, we shall write $\|\cdot\|_{\Sigma}$ instead of $\|\cdot\|_{\M_\Sigma}=\|\cdot\|_{2, \mathcal P_\Sigma}$ and $\|\cdot\|_{B, \Sigma}$ instead of $\|\cdot\|_{\M^B_\Sigma}$. Also, if $\Sigma$ is the trivial partition of $\s(A)$ into singletons, we may write simply $\M^B$ instead of $\M^B_\Sigma$.

First, given $X\in\mathfrak S_2(H, \mathcal P_\Sigma)$, we  define a sequence $(\alpha_n(X))_{n\in\ZZ}$   by
\begeq\label{alpha}
\alpha_n(X)=\|X\|_\Sigma^{-\frac{1}{2}}\max\left\{\left(\sum\limits_{\substack{|k|\ge  |n| \\
k\in\mathbb{Z}}}\|P_{\s_k}X\|_\Sigma^2\right)^\frac{1}{4}, \left(\sum\limits_{\substack{|k|\ge  |n| \\
k\in\mathbb{Z}}}\|XP_{\s_k}\|_\Sigma^2\right)^\frac{1}{4}\right\}. 
\eq
It is easy to check that the above sequence has the following properties:
\ben
\item $\alpha_n(X)=\alpha_{-n}(X)$, $n\in\mathbb{Z}$.

\item $\lim\limits_{|n|\to\infty}\alpha_n(X)=0$, $n\in\mathbb{Z}$.

\item $\alpha_n(X)\leq  1$ for all $n\in\mathbb{Z}$.

\item $\alpha_n(X)\geq  \alpha_{n+1}(X)$, $n\geq  0$.
\een
\noindent In addition, if we  let $P^\Sigma_{(m)} = \sum_{|n|\le m} P_{\s_n}$ and assume that $P^\Sigma_{(m)}XP^\Sigma_{(m)}\ne X$ for all $m\in\mathbb{Z}_+$,   we have
\ben
\item[5.] $\alpha_n(X)\ne 0$ for all $n\in\mathbb{Z}$, and
\begeq\label{weight}
\sum\limits_{n\in\mathbb{Z}}\frac{\|XP_n\|_\Sigma^2+\|P_nX\|_\Sigma^2}{(\alpha_n(X))^2} <\infty.
\eq
\een 

We remark that the sequence $(\a_n(X))$ characterizes the decay of the entries of $X^\Sigma$ along its rows and columns. In view of \eqref{weight}, one may conclude that any $X\in\M_\Sigma$ also belongs to a ``weighted Hilbert-Schmidt space'' with a weight that depends on $X$. This is a manifestation of the fact that for any convergent series there is another one albeit with a slower rate of convergence.  If $X = B$, we shall write simply $\a_n$ instead of $\a_n(B)$.

We also remark that if the perturbation $B$ is such that  $P^\Sigma_{(m)}BP^\Sigma_{(m)} = B$ for some $m\in\ZZ_+$, then there is no need for the space $\M^B_\Sigma$ in the spectral analysis of $A-B$. Therefore, in this section, we assume that \eqref{weight} holds for $X = B$.

Next, we introduce the function $f = f_B: \sigma(A)\to\mathbb{R}_+$ given by
$$
f(\lambda)=\sum_{n\in\ZZ}\alpha_n\chi_{\s_n}(\l), \ \l\in\s(A),
$$
where 
$\chi_E$ denotes the characteristic function of a set $E$.  Using the functional calculus for unbounded normal operators \cite{DS88II}, we get that the operator
$$
f(A)=\sum_{n\in\mathbb{Z}}\alpha_nP_{\s_n}
$$
belongs to $B({H})$ and $\|f(A)\|\leq  \max\limits_{n\in\mathbb{Z}}|\alpha_n|=1$.

We let $\M^B_\Sigma$ be the set of all operators  $X\in \M_\Sigma$ such that there 
exist operators $X_l, X_r\in\M_\Sigma$ satisfying
\begeq\label{xlxr}
X=X_lf(A) \quad\mbox{and}\quad X=f(A)X_r.
\eq
Observe that the $\Sigma$-matrix of the operator $f(A)$ is diagonal and the assumption $B\neq  P_{(m)}BP_{(m)}$ for all $m\in\ZZ_+$ implies that it  has no zeros on the main diagonal (see the last property of the sequence $\a_n$). Therefore, given $X \in\M$, the operators $X_l$ and $X_r$ are uniquely determined by \eqref{xlxr}. Moreover, $\M^B_\Sigma$ is a Banach space with the norm
$\|X\|_{B,\Sigma} = \max\{\|X_l\|_\Sigma, \|X_r\|_\Sigma\}$ and 
\begeq\label{Membed}
\|X\|_\Sigma = \|X_lf(A)\|_\Sigma = \|f(A)X_r\|_\Sigma \le \|X\|_{B,\Sigma}.
\eq
In fact, $\M^B_\Sigma$ is a Banach algebra and
\begeq\label{mbs}
\|XY\|_{B,\Sigma} \le \min\left\{\|X\|_\Sigma\|Y\|_{B,\Sigma}, \|Y\|_\Sigma\|X\|_{B,\Sigma}\right\},\
 X,Y\in \M^B_\Sigma.
\eq 

From \eqref{weight}, we also deduce that $B\in\M^B_\Sigma$ with 
$$B_l = \sum_{n\in\ZZ}\frac1{\a_n} BP_{\s_n} \mbox{ and } B_r =  \sum_{n\in\ZZ}\frac1{\a_n} P_{\s_n}B.$$

 
 Next, we will show that for any partition $\wS$ that is coarser than $\Sigma$, we have that 
$(\M^B_\Sigma, J_\wS, \Gamma_{\wS})$ is an admissible triplet. Moreover, we will see  that for any $\epsilon > 0$ there exists a partition $\wS$ such that the constant $\gamma = \g_\wS$ in Property \ref{adprop4} of Definition \ref{baskdef8} satisfies  
$\g_\wS < \epsilon.$
This allows one to use Theorem \ref{baskth6} with \eqref{bask12} automatically satisfied.

\bp\label{trM}
Assume that the operator $A$ satisfies Assumption \ref{assa},  $\Sigma$ is a partition of the spectrum $\s(A)$, and $B\in\M_\Sigma$. For any partition $\wS$ that is coarser than $\Sigma$, we have that $(\M^B_\Sigma, J_\wS, \Gamma_{\wS})$ is an admissible triplet.
Additionally, for $m\in\ZZ_+$, let  $\widetilde\s_m = \bigcup_{|n|\le m} \s_n$, 
\begeq\label{sigM}
\Sigma_m = \left\{\widetilde\s_m\right\}\cup\left\{\s_n\in\Sigma: |n|> m\right\},
\eq
and $\g_m$ be the constant   in Property \ref{adprop4} of Definition \ref{baskdef8} for 
$(\M^B_\Sigma, J_{\Sigma_m}, \Gamma_{\Sigma_m})$. Then one can choose $\g_m$ so that
\begeq\label{limgam}
\lim\limits_{m\to \infty} 
\g_m =0.
\eq
\ep

\bpf
The first and last properties of the admissible triplet in Definition \ref{baskdef8} follow from \eqref{Membed} and Assumption \ref{assa}(\ref{assa5}).
Next, observe that
\begeq\label{Jinv}
J_\wS(X_lf(A))=(J_\wS X_l)f(A), \quad J_\wS(f(A)X_r)=f(A)(J_\wS X_r),
\eq
\begeq\label{Ginv}
\Gamma_\wS(X_lf(A))=(\Gamma_\wS X_l)f(A), \quad\mbox{and}\quad \Gamma_\wS(f(A)X_r)=f(A)(\Gamma_\wS X_r),
\eq
for $X_r$, $X_l\in\M_\Sigma$. It follows that the space $\M^B_\Sigma$ is invariant for $J_\wS$ and $\Gamma_\wS$ and, therefore,  the restrictions of the transforms to $\M^B_\Sigma$ are well defined. The 
remaining properties of the admissible triplet easily follow. 

To ensure \eqref{limgam}, we need to obtain better estimates for $\gamma_m$ than follow directly from \eqref{mbs} and \eqref{Ginv}. To do so, we shall make use of two more sequences:
$(\alpha_n')$ and $(\widetilde{\alpha}_n)$, $n\in\NN$. The first of them is defined by
$$
\alpha_{n+1}'=\sup
\left\{\alpha_\ell\max\{d_{j\ell}, d_{\ell j}\}: \ell, j\in\mathbb{Z}, |\ell|\le n, |j|> n\right\}, \quad n\in\mathbb{Z}_+,
$$
where $$d_{j\ell}= 
\left(\sup_{\l_\ell\in\s_\ell} \sum_{\l_j\in\s_j} |\lambda_j-\l_\ell|^{-2}\right)^{\frac12},\quad
\ \ell, j\in\mathbb{Z},\ j\neq\ell.$$ 
The second sequence is given by
\begeq\label{altil}
\widetilde{\alpha}_n=\sqrt \eta\alpha_n+\alpha_n', \quad n\in\NN,
\eq
where $\eta$ is defined by \eqref{specAcond}. 
 Observe that since $$\lim_{n\to \infty} d_{j,j+n} =\lim_{n\to \infty} d_{j+n,j} = 0,\quad j\in\ZZ,$$ we have
 \begeq\label{limn} 
 \lim\limits_{n\to\infty}\alpha_n'= \lim\limits_{n\to\infty}\widetilde\alpha_n=0.
 \eq
 
 We will also need the estmates
 \begeq\label{gxf}
 \left\|(\Gamma_{\Sigma_m}(Xf(A))^\Sigma_{j\ell} \right\|_{B(H)} \le
\a_\ell d_{j\ell}\|X^\Sigma_{j\ell}\|_{B(H)},\quad X\in\M_\Sigma,
\eq
and
\begeq\label{gfx}
 \left\|(\Gamma_{\Sigma_m}(f(A)X)^\Sigma_{j\ell} \right\|_{B(H)} \le
\a_j d_{j\ell}\|X^\Sigma_{j\ell}\|_{B(H)},\quad X\in\M_\Sigma,
\eq
which follow the same way as in the proof of the bound for $\|\Gamma\|_{B(B(H))}$.

Firstly, we obtain a bound for $\|\Gamma_{\Sigma_m}\|_{\M^B_\Sigma\to B(H)}$.

Let $Q_{\widetilde\s_m}=I-P_{\widetilde\s_m}$, $m\in\mathbb{Z}_+$. Then,
for $X\in\M_\Sigma$, we use 
\eqref{Gs} and \eqref{gxf} to get
\begeq\label{gga1}
\bs
\|\Gamma_{\Sigma_m}(Xf(A))\|_\Sigma & = \|\Gamma_{\Sigma_m}(Xf(A)Q_{\widetilde\s_m}) + \Gamma_{\Sigma_m}(Q_{\widetilde\s_m}Xf(A)P_{\widetilde\s_m})\|_\Sigma\\ 
&\le(\sqrt \eta\a_{m+1} +\a'_{m+1})\|X\|_\Sigma = \widetilde\a_{m+1}\|X\|_\Sigma.
\end{split}
\eq
Similarly, we have
\begeq\label{gga2}
\bs
\|\Gamma_{\Sigma_m}(Xf(A))\|_\Sigma & = \|\Gamma_{\Sigma_m}(Q_{\widetilde\s_m}f(A)X) + \Gamma_{\Sigma_m}(P_{\widetilde\s_m}Xf(A)Q_{\widetilde\s_m})\|_\Sigma\\ 
&\le(\sqrt \eta\a_{m+1} +\a'_{m+1})\|X\|_\Sigma = \widetilde\a_{m+1}\|X\|_\Sigma.
\end{split}
\eq
From \eqref{gga1} nad \eqref{gga2}, it follows that
\begeq\label{gm1}
\|\Gamma_{\Sigma_m}\|_{\M^B_\Sigma\to B(H)}\le \|\Gamma_{\Sigma_m}\|_{\M^B_\Sigma\to \M_\Sigma} \le \widetilde\a_{m+1}.
\eq

Secondly, we need  estimates for $\|X\Gamma_{\Sigma_m} Y\|_{B,\Sigma}$ and $\|(\Gamma_{\Sigma_m} X)Y\|_{B,\Sigma}$. However, from
\eqref{mbs} and \eqref{gm1}, we immediately get
\begeq\label{gm2}
\max\{\|X\Gamma_{\Sigma_m} Y\|_{B,\Sigma}, \|(\Gamma_{\Sigma_m} X)Y\|_{B,\Sigma}\} \le \widetilde\a_{m+1}\|X\|_{B,\Sigma}\|Y\|_{B,\Sigma}.
\eq
Thus, we can take $\g_m =  \widetilde\a_{m+1}$ and \eqref{limgam} follows from \eqref{limn}.
\epf

\section{Main results}\label{main}

In this section, we collect various versions of the similarity results for the operator $A-B$ that follow from the more abstract theorems of the previous section. To simplify the exposition, we collect the generic conclusions of the theorems below in two statements. In the first statement,  we denote by $\M$ some space of admissible perturbations and by $\Sigma = \{\s_n, n\in\mathcal I\}$  -- some partition of the spectrum $\s(A)$. In the theorems, the concrete $\M$ and $\Sigma$ will be provided.
In all of the theorems, the operator $A$ is assumed to satisfy Assumption \ref{assa}.

\bcon\label{con}
There exist operators $U, V\in\M$ such that the operator $A-B$ is similar to the operator $A-V$, and the similarity transform is given by the operator $I+U$:
\begeq\label{msimeq}
(A-B)(I+U) = (I+U)(A-V).
\eq
Moreover, the $\Sigma$-matrix of the operator $V$ is diagonal, and the spectral projections $ P^\prime_{\s_n^\prime} = (I+U)P_{\s_n}(I+U)^{-1}$ of the operator $A-B$ form a resolution of the identity. 
\econ

The second conclusion provides a rough asymptotic estimate of the spectrum $\s(A-B)$. Given a weight function $w: \ZZ\to [1,\infty)$, we write $\ell^2_w = \ell^2_w(\ZZ)$ for the weighted $\ell^2$ space of sequences $x = (x_n)\in \ell^2$ such that $\sum_{n\in\ZZ} |x_n|^2w(n) < \infty$.

\bcon\label{con2} We have
\begeq\label{spec1}
\s(A-B) = \bigcup_{n\in\ZZ} \{\l_n - b_n\}
\eq
for some sequence $b = (b_n)\in\ell^2_w(\ZZ)$.
\econ

Particular weights $w$ will also be specified in the theorems below.


\subsection{Results for  perturbations in $\mathfrak S_2(H, \mathcal P_\Sigma)$.}

In this subsection, we deal with perturbations $B\in\mathfrak S_2(H, \mathcal P_\Sigma)$ so that we can avoid using Theorem \ref{baskth8}. We continue to use the notation of Section \ref{method}.

We begin with small perturbations $B$ in order to make
Theorem \ref{baskth6} directly applicable for $A-B$. 

\bt\label{mt1}
Assume that the operator $A$   satisfies Assumption \ref{assa} and $B\in\mathfrak S_2(H)$ with
$\|B\|_2 < \frac\delta4$, where $\delta$ is given by \eqref{speccond}. Then Conclusion \ref{con}   holds with $\M = \mathfrak S_2(H)$,  $U = \Gamma X_*$, and $V = JX_*$, where $X_*\in\M$ is the limit of simple iterations: $X_0=0$, $X_1=\Phi(X_0) = B$, etc.~with $\Phi$ given by \eqref{bask13}. Moreover, Conclusion  \ref{con2}  holds with the trivial weight $w \equiv 1$.
\et

\bpf
Conclusion \ref{con} follows immediately from Theorem \ref{baskth6} and Lemmas \ref{bask100} and \ref{JG1}. To obtain Conclusion \ref{con2} we also need Lemma \ref{lemspecsum}.
\epf

\bt\label{mt2}
Assume that the operator $A$   satisfies Assumption \ref{assa},  $B\in\mathfrak S_2(H, \mathcal P_\Sigma)$ for some partition $\Sigma$ of $\s(A)$, and
$\|B\|_\Sigma < \frac1{4\sqrt \eta}$, where $\eta$ is given by \eqref{specAcond}. Then Conclusion \ref{con} holds with $\M = \mathfrak S_2(H, \mathcal P_\Sigma)$,  $U = \Gamma_\Sigma X_*$, and $V = J_\Sigma X_*$, where $X_*\in\M$ is the limit of simple iterations: $X_0=0$, $X_1=\Phi(X_0) = B$, etc.~with $\Phi$ given by \eqref{bask13}.
\et

\bpf
Follows immediately from Theorem \ref{baskth6} and Lemmas \ref{bask100} and \ref{JG2}.
\epf

In view of Remark \ref{rem3}, we have the following slight modifications of the above two results.

\bt\label{mt11}
Assume that the operator $A$   satisfies Assumption \ref{assa}, and $B\in\mathfrak S_2(H)$ with $JB =0$ and
$\|B\|_2 < \frac\delta3$, where $\delta$ is given by \eqref{speccond}. Then Conclusion \ref{con}  holds with $\M = \mathfrak S_2(H)$,  $U = \Gamma X_*$, and $V = JX_* = J(B\Gamma X_*)$, where $X_*\in\M$ is the limit of simple iterations: $X_0=0$, $X_1=\Phi(X_0) = B$, etc.~with $\Phi$ given by 
$\Phi(X) = B\Gamma X-(\Gamma X)J(B\Gamma X)+B,$ $X\in\M$. Moreover, Conclusion \ref{con2} holds for some sequence $b\in\ell^1(\ZZ)$. 
\et

\bpf
The equation $JX_* = J(B\Gamma X_*)$ follows from \eqref{JXs}. The fact that $b\in\ell^1$ follows since $B\Gamma X_* \in \mathfrak S_1(H)$ as a product of two operators in $\mathfrak S_2(H)$.
\epf

\bt\label{mt21}
Assume that the operator $A$   satisfies Assumption \ref{assa},  $B\in\mathfrak S_2(H, \mathcal P_\Sigma)$ for some partition $\Sigma$ of $\s(A)$, $J_\Sigma B = 0$, and
$\|B\|_\Sigma < \frac1{3\sqrt \eta}$, where $\eta$ is given by \eqref{specAcond}. Then Conclusion \ref{con} holds with $\M = \mathfrak S_2(H, \mathcal P_\Sigma)$,  $U = \Gamma_\Sigma X_*$, and $V = J_\Sigma X_* = J_\Sigma (B\Gamma_\Sigma X_*)
$, where $X_*\in\M$ is the limit of simple iterations: $X_0=0$, $X_1=\Phi(X_0) = B$, etc.~with $\Phi$ given by 
$\Phi(X) = B\Gamma_\Sigma X-(\Gamma_\Sigma X)J_\Sigma(B\Gamma_\Sigma X)+B,$ $X\in\M$.
\et

Next, we use Proposition \ref{trM} to strengthen the first two of the above four theorems by removing the condition on the norm of the perturbation.


\bt\label{mt12}
Assume that the operator $A$   satisfies Assumption \ref{assa} and  $B\in\mathfrak S_2(H, \mathcal P_\Sigma)$ for some partition $\Sigma$ of $\s(A)$. Then there exists $m\in\ZZ_+$ such that Conclusion \ref{con}   holds with $\M = \M^B_{\Sigma_m}$,  $U = \Gamma_{\Sigma_m} X_*$, and $V = J_{\Sigma_m}X_*$, where $\Sigma_m$ is given by \eqref{sigM}, and $X_*\in\M$ is the limit of simple iterations: $X_0=0$, $X_1=\Phi(X_0) = B$, etc.~with $\Phi$ given by \eqref{bask13}. Moreover, if $\Sigma$ is the trivial partition of $\s(A)$ into singletons, then Conclusion  \ref{con2}  holds  
with the weight $w$ given by $w(k) = (\a_k(B))^{-2}$, where $\a$ is defined by \eqref{alpha}.
\et

\bpf
Conclusion \ref{con} follows immediately from Theorem \ref{baskth6} and Proposition \ref{trM}. 
To obtain Conclusion \ref{con2} we again need Lemma \ref{lemspecsum} and, also, the fact that  an operator with a diagonal matrix is in $\M^B$ if and only if the sequence of elements on the diagonal is in $l^2_w$ with $w(k) = (\a_k(B))^{-2}$. In case of the trivial partition $\Sigma$, the operator $J_{\Sigma_m}X_*\in \M^B_{\Sigma_m}$  is a finite-rank perturbation of the diagonal operator $JX_*\in\M^B$, and the result follows.
\epf

\subsection{Results using the preliminary similarity transform.} 

In this subsection, we deal with perturbations $B$ which require the use of both the preliminary transform (Theorem \ref{baskth8}) and the weighted Hilbert-Schmidt space introduced in Subsection \ref{wspace}. 

To use Theorem \ref{baskth8}, we need to verify that Assumption \ref{baskpred1} holds. Here we shall take $\M = \mathfrak S_2(H)$. In lieu of the transforms in the admissible triplet, we shall use the coarser transforms $J_m$ and $\Gamma_m$ that correspond to the partition $\mathcal P_{(m)}$ in Example \ref{pex}, $m\in\ZZ_+$. The transforms are defined by \eqref{tJ2s} and \eqref{Gs} with the choice of $\Sigma = \Sigma_{m}$ that generates $\mathcal P_{(m)}$. 

\bl
Assume that the operator $A$   satisfies Assumption \ref{assa} and 
the operator $B$ satisfies Assumptions \ref{asB} and \ref{asBb}. Then there exists $m\in\ZZ_+$ such that Assumption \ref{baskpred1} holds with $\M = \mathfrak S_2(H)$, $J = J_m$ and $\Gamma = \Gamma_m$.
\el

\bpf
Observe that \eqref{cGB} is equivalent to $\Gamma_0 B\in\M$, \eqref{cBGB} -- to $B\Gamma_0 B\in\M$, and Assumption \ref{asBb} -- to $J_0 B\in B(H)$.\footnote{Assumption \ref{asBc}, which we did not impose in this lemma, is equivalent to $J_0 B\in\M$.} Observe also that for each $m\in\ZZ_+$, we have that $J_mB$ and $\Gamma_m B$ are finite-rank perturbations of $J_0B$ and $\Gamma_0 B$, respectively. Therefore,  Assumption \ref{asB}\eqref{assb2} and Assumption \ref{asBb} yield Assumption \ref{baskpred1}\eqref{mprop}. Moreover, Assumption \ref{asB}\eqref{assb3} yields  Assumption \ref{baskpred1}\eqref{propas5}. Properties \ref{propas2} and \ref{propas4} of Assumption \ref{baskpred1} are obtained the same way as in the proofs of Lemma \ref{bask100} and \ref{JG2}. It remains to show that there exists $m\in\ZZ_+$ such that $\|\Gamma_m B\|_2 < 1$. The latter, however, follows immediately from the definition of $\Gamma_m$ since $\Gamma_0 B\in\M$. Indeed, we have
\begeq\label{Gm0}
\lim_{m\to \infty} \|\Gamma_m B\|_2 =0,
\eq
and the lemma is proved.
\epf

Immediately from the above lemma and Theorem \ref{baskth8}, we get the following result.

\bt\label{1pre}
Assume that the operator $A$   satisfies Assumption \ref{assa} and the operator $B$ satisfies Assumptions \ref{asB} and \ref{asBb}. Then there exists $m\in\ZZ_+$ such that
the operator $A-B$ is similar to   $A-J_mB-B_0$, where
$B_0=(I+\Gamma_m B)^{-1}(B\Gamma_m B-(\Gamma_m B)J_m B)$. 
The similarity transform is given by $I+\Gamma_m B$ so that
\begeq\label{abs2}
(A-B)(I+\Gamma_m B)=(I+\Gamma_m B)(A-J_mB-B_0).
\eq
We also have $B_0 \in\mathfrak S_2(H)$ and $B_0- B\Gamma_0 B\in\mathfrak S_1(H)$. Moreover, if we additionally impose Assumption \ref{asBc}, then $J_0 B$ and  $J_m B$ are in $\mathfrak S_2(H)$.
\et

\bpf
The only statement in the theorem that we have not discussed before is $B_0- B\Gamma_0 B\in\mathfrak S_1(H)$. This, however, follows immediately from the Neumann series representation of 
$(I+\Gamma_m B)^{-1}$ and the fact that $B\Gamma_0 B$ is a finite-rank perturbation of $B\Gamma_m B$.
\epf

Next, we would like to apply Theorem \ref{baskth6} to the operator $A-Q$, where $Q = J_mB + B_0$  was constructed in the above theorem.
We have two cases. If we impose Assumption \ref{asBc} on $B$, then
 $Q \in \M = \mathfrak S_2(H)$, and the way is fairly straightforward in view of Proposition \ref{trM}. 
 Alternatively, we will have to run the scheme for the operator $A - J_m B$ instead of $A$ and that would require more work and additional assumptions on $B$.

\bt\label{mt3}
Suppose that the operator $A$   satisfies Assumption \ref{assa} and  the operator $B$ satisfies Assumptions \ref{asB} and \ref{asBc}. Then there exist $m\le k\in\ZZ_+$ such that Conclusion \ref{con}
holds with
\ben
\item $ U = \Gamma_m B +\Gamma_k X_*+(\Gamma_m B)(\Gamma_k X_*)$;
\item $V = J_k X_* = J_k Q + J_k(Q\Gamma_k X_*)  = J_m B + J_k(B_0(I+\Gamma_k X_*))$;
\een
where
\ben
\item $Q = J_m B + B_0$; 
\item $ B_0 = (I+\Gamma_m B)^{-1}(B\Gamma_m B - (\Gamma_m B)J_m B)$;
\item $X_*\in\M$ is the limit of simple iterations: $X_0=0$, $X_1=\Phi(X_0) = Q$, etc.~with $\Phi$ given by 
$\Phi(X) = Q\Gamma_k X-\Gamma_k X(J_kQ) - (\Gamma_k X)J_k(Q\Gamma_k X)+Q,$ $X\in\M$. 
\een
Moreover, we have 
\begeq\label{formV}
V = J_0 B + J_0(B\Gamma_0 B)+C
\eq
 for some $C\in \mathfrak S_1(H)$, and Conclusion \ref{con2} holds with $w(\ell) = (\a_\ell(Q))^{-2}$, where $\a$ is defined by \eqref{alpha}. 
\et

\bpf
The result follows immediately from Proposition \ref{trM} (with $Q$ in place of $B$) and Theorems \ref{baskth6} and \ref{1pre}. The formulas for $V$ are easily obtained by direct computation using \eqref{JXs} and the 
definitions of $\Gamma_k$ and $J_k$, $k\in\ZZ_+$. Conclusion \ref{con2} follows as in the proof of Theorem \ref{mt12}.
\epf


\bc\label{cor1d}
With the notation of Theorem \ref{mt3}, assume in addition that all spectral projections $P_n$, $n\in\ZZ$, have rank $1$. Then there exist sequences $p,q \in \ell^2(\ZZ)$ and $c\in\ell^1(\ZZ)$ such that $(J_0B)P_n = P_n(J_0B)P_n = p_n P_n$ and $(J_0(B\Gamma_0 B))P_n = P_n(J_0(B\Gamma_0 B))P_n = q_n P_n$, $n\in\ZZ$, and the sequence $b$ in Conclusion \ref{con2} satisfies $b = p+q+c$.
\ec

\bpf
The result follows immediately from \eqref{formV}.
\epf

\brem
For a specific perturbation $B$ the sequences $p$ and $q$ in the above corollary are explicitly computable.
\erem

Let us now discuss the case when $J_m B \in B(H)\setminus\mathfrak S_2(H)$. Here, we would like to treat $A -J_mB$ as an unperturbed operator. This necessitates the following additional assumption.

\begin{ass}\label{abba}
For any $N\in \mathbb N$ there exists $m > N$ such that the operator $A - J_m B$ satisfies Assumption \ref{assa}.
\end{ass}

Under Assumption \ref{abba}, we have the resolution of the identity $\wcP = \{\wP_n, n\in\ZZ\}$ that consists of the spectral projections of the operator $A-J_m B$. This partition and its coarser versions allow us to define two families of transforms $\widetilde J_k$ and $\widetilde\Gamma_k$, $k\in\ZZ$, for the operator $A - J_m B$ the same way $J_k$ and $\Gamma_k$, $k\in\ZZ$, were defined for $A$.

Since $B_0 \in\mathfrak S_2(H)$, where $B_0$ was constructed in Theorem \ref{1pre}, we have that Proposition \ref{trM} applies for $A-J_mB$ in place of $A$ and $B_0$ in place of $B$. This yields the admissible triplets $(\M^{B_0}_{\wS_0}, \widetilde J_k, \widetilde \Gamma_k)$, $k\in\ZZ_+$, where
$\wS_0$ is the partition of $\s(A-J_m B)$ into singletons. Applying Theorem \ref{baskth6} for the free operator $A-J_m B$, perturbation $B_0$, and a triplet $(\M^{B_0}_{\wS_0}, \widetilde J_k, \widetilde \Gamma_k)$ with a sufficiently large $k\in\ZZ_+$, we deduce the following result.



\bt\label{mt4}
Suppose that the operator $A$  and $B$ satisfy Assumptions \ref{assa}, \ref{asB},  \ref{asBb}, and \ref{abba}. Then there exist $m, k\in\ZZ_+$ such that Conclusion \ref{con}
holds with $\M = B(H)$ and
\ben
\item $ U = \Gamma_m B +\widetilde\Gamma_{k} X_*+(\Gamma_m B)(\widetilde\Gamma_{k} X_*) \in \mathfrak S_2(H)$;
\item $V = J_m B + \widetilde J_{k} X_*  = J_m B + \widetilde J_{k}(B_0(I+\widetilde\Gamma_{k} X_*))$;
\een
where
\ben
\item $ B_0 = (I+\Gamma_m B)^{-1}(B\Gamma_m B - (\Gamma_m B)J_m B) \in\mathfrak S_2(H)$;
\item $X_*\in\mathfrak S_2(H)$ is the limit of simple iterations: $X_0=0$, $X_1=\Phi(X_0) = B_0$, etc.~with $\Phi$ given by 
$$\Phi(X) = B_0\widetilde\Gamma_{k} X-\widetilde\Gamma_{k} X(\widetilde J_{k}B_0) - (\widetilde\Gamma_{k} X)\widetilde J_{k}(B_0\widetilde\Gamma_{k} X)+B_0, \ X\in\mathfrak S_2(H).$$ 
\een
Moreover, we have 
\begeq\label{formV1}
V = J_0 B + \widetilde J_{0}(B\Gamma_0 B)+C
\eq
 for some $C\in \mathfrak S_1(H)$. 
\et


\brem\label{badest}
Conclusion \ref{con2} does not have to hold in the above theorem. In many examples \cite{BKU19}, however, Formula \eqref{formV1} can still be used effectively to estimate the spectrum $\s(A-B)$.
\erem

\brem
In some cases, a version of Theorem \ref{mt4} would remain true even if $J_m B \in \mathfrak L_A(H)\setminus B(H)$. One only needs to ensure that $(\Gamma_m B)(J_m B)\in\mathfrak S_2(H)$.
\erem

\subsection{Equiconvergence of  spectral decompositions.}

In this subsection, we explore the consequences of Conclusion \ref{con} for approximation of the spectral projections $P^\prime_{\s^\prime}$ of the operator $A-B$, where  $\s^\prime$ is a spectral component of $\s(A-B)$. Recall that, if the conclusion holds, we have $ P^\prime_{\s^\prime} = (I+U)P_{\s}(I+U)^{-1}$, where $\s$ is the corresponding component of $\s(A)$.  
It follows that
\begeq\label{prodif}
 P^\prime_{\s^\prime} - P_{\s}= (UP_{\s}- P_{\s}U)(I+U)^{-1} \in \mathfrak S_2(H),
\eq
whenever $U\in\mathfrak S_2(H)$. 

As was the case in most of the major results above, we assume now that $U \in \M^Q_\Sigma$ for some fixed partition $\Sigma = \{\s_n: n\in\Omega\}$ of $\s(A)$ and $Q\in \mathfrak S_2(H, \mathcal P_\Sigma)$. Further estimates shall be obtained in terms of the sequence $\a_n(Q)$ defined by \eqref{alpha} . Furthermore, given a spectral component $\s$ of $\s(A)$, we let
\[
\a_n^\s = \a_n^\s(Q) = \max\{\a_n(Q): \s_n\cap\s\neq\emptyset\}.
\]

\bl\label{proin}
With the above notation, we have
\[
\max\{\|UP_{\s}\|_\Sigma, \|P_{\s}U\|_\Sigma \} \le \a_n^\s \|U\|_{Q,\Sigma}.
\]
\el
\bpf
We have
\[
\bs
\|UP_{\s}\|_\Sigma^2 &= \left\| \sum_{n\in\ZZ}\a_n(Q)U_lP_{\s_n}P_{\s} \right\|^2_\Sigma =
\sum_{n\in\ZZ}\a^2_n(Q)\left\| U_lP_{\s_n}P_{\s} \right\|^2_\Sigma \\
& \le (\a_n^\s)^2 \sum_{n\in\ZZ}\left\| U_lP_{\s_n}P_{\s} \right\|^2_\Sigma \le (\a_n^\s)^2 \|U\|_{Q,\Sigma}^2.
\end{split}
\]
The second inequality is obtained in a similar fashion.
\epf

\bt\label{tpro}
Assume that $\|U\|_\Sigma < 1$. With the above notation, we have
\[
\|P^\prime_{\s^\prime} - P_{\s}\|_\Sigma  \le \frac{2\|U\|_{Q,\Sigma}}{1-\|U\|_\Sigma}\a_n^\s.
\]
\et

\bpf
The  estimate  follows immediately from \eqref{prodif}, Lemma \ref{proin}, and the inequality
\[
\|(I+U)^{-1}-I\|_\Sigma \le \frac{\|U\|_\Sigma}{1-\|U\|_\Sigma},
\]
which is implied by the Neumann series representation of $(I+U)^{-1}$.
\epf

We conclude this subsection with the following result on equiconvergence of the spectral decompositions. 
\bt
Assume that Conclusion \ref{con} holds with $U \in \M^Q_\Sigma$ for some fixed partition $\Sigma = \{\s_n: n\in\Omega\}$ of $\s(A)$ and $Q\in \mathfrak S_2(H, \mathcal P_\Sigma)$. Let $\s_{(n)} = \s(A) \setminus \bigcup_{|m|< n} \s_m$ and $\s^\prime_{(n)}$ be the corresponding spectral component of $\s(A-B)$, $n\in\ZZ$. 
Assume also that $\|U\|_\Sigma < 1$ or $\Sigma$ is the trivial partition of $\s(A)$ into singletons. Then
\[
\lim_{n\to\infty} \left\|P^\prime_{\s(A-B)\setminus\s^\prime_{(n)}} - P_{\s(A)\setminus\s_{(n)}}\right\|_\Sigma =\lim_{n\to\infty} \left\|P^\prime_{\s^\prime_{(n)}} - P_{\s_{(n)}}\right\|_\Sigma = 0.
\]
\et
\bpf
The result immediately follows from Theorem \ref{tpro} because $\a_n^{\s_{(n)}} = \a_n(Q) \to 0$ as $n\to\infty$. In the case of the trivial partition $\Sigma$, the condition $\|U\|_\Sigma < 1$ is unnecessary because one can use \eqref{prodif} and the fact that $\|MN\|_2 \le \|M\|_2\|N\|$ for any pair of operators $M\in\mathfrak S_2(H)$ and $N\in B(H)$.
\epf
\subsection{Spectral splitting.}

In this subsection, we run the method  of similar operators in the case when the spectrum $\s(A)$ is partitioned into just two parts: $\Sigma = \{\s_1, \s_2\}$ with $\s_1 = \{\l_k\}$ and $\s_2 = \s(A)\setminus \{\l_k\}$.  
Therefore, the partition $\mathcal P_\Sigma$ is one of the partitions that appeared in Example \ref{pex}: $\mathcal P_\Sigma = \mathcal P_k = \{P_k, Q_k\}$, where $Q_k = I - P_k$. 
The result is a special case of the spectral splitting, the theory of which have been developed extensively in \cite{B86}. The perturbation $B$ here is assumed to be in $B(H)$ and sufficiently small.

For   the second result in this subsection, we will also assume that
 $\l = \l_k$ is a simple eigenvalue of $A$, i.e.~the corresponding eigenspace is spanned by the single eigenvector $e = e_k$. The goal then is to obtain good estimates for the eigenvalue $ \l^\prime =\l_k^\prime$ of the operator $A-B$ and the corresponding eigenvector $ e^\prime = e^\prime_k$. To simplify the notation we shall also use $P = P_k$ and $Q=Q_k$ whenever no ambiguity may arise. Thus,
for $\Sigma$-matrices in $\mathfrak L_A(H, \mathcal P_\Sigma)$ we shall write:
\[
X = \left( \begin{array}{cc}
X_{11} & X_{12} \\
X_{21}& X_{22}
\end{array} 
\right) = \left( \begin{array}{cc}
PXP & PXQ \\
QXP& QXQ
\end{array} 
\right).
\]
Formula \eqref{tJ2s}, defining the transform $J = J_\Sigma$ 
reduces, in this case, to
\[
JX = \left( \begin{array}{cc}
X_{11} & 0 \\
0& X_{22}
\end{array} 
\right),\quad X\in \mathfrak L_A(H).
\]
The formula \eqref{Gs} for the transform $\Gamma = \Gamma_\Sigma$ can also be written in a different way. There exists an operator $S\in B(H)$ defined via
\begeq\label{opS}
PS = SP = 0,\quad S(\l I-A)= (\l I-A)S = Q. 
\eq
A straightforward computation shows that the transform $\Gamma= \Gamma_\Sigma$ defined by \eqref{Gs} satisfies
\begeq\label{Gsp}
\Gamma X = PXS-SXP,\quad X\in B(H).
\eq

Applying Lemma \ref{JG2}, we see that the  triplet 
$(\M, J, \Gamma)$ with $\M = B(H)$ 
 is admissible. 
 Observe also that from \eqref{Gsp} we immediately have
\begeq\label{nG}
\|\Gamma\|_{B(\M)} \le  s\sqrt2, \quad \mbox{where}\quad s =\|S\| = \sup\limits_{j\in\ZZ\setminus\{k\}} |\l_k-\l_j|^{-1} = \delta^{-1}.
\eq

Applying Theorem \ref{baskth6}  to $(\M, J, \Gamma)$ we get the following result, a stronger version of which can be found in \cite{B86}.

\bt\label{mt5b}
Suppose that the operator $A$   satisfies Assumption \ref{assa}, and  the operator $B\in \M = B(H)$ is such that $\|B\|_\M < \frac1{4s\sqrt2}$. Then Conclusion \ref{con} holds with $U = \Gamma X_*$ and $V = JX_* = X_{11}^*+ X_{22}^*$, where $X_* =(X_{ij}^*)$, $i,j =1,2$, is the solution of the following non-linear system of  equations:
\begeq\label{sXs}
\left\{\begin{array}{l}
X_{11} = - B_{12}SX_{21}+B_{11}\\
X_{21} = - B_{22}SX_{21}+SX_{21}B_{11}-SX_{21}B_{12}SX_{21}+B_{21}\\
X_{12} = B_{11}X_{12}S - X_{12}SB_{22} - X_{12}SB_{21}X_{12}S+B_{12}\\
X_{22} = B_{21}X_{12}S+B_{22}.
\end{array}
\right.
\eq
\et

\bpf
One easily verifies by direct computation that in this case the $\Sigma$-matrix  $X_*$ solves \eqref{sXs} if and only if it is a fixed point of the function $\Phi$ given by \eqref{bask13}.
\epf

\brem
In \cite[Theorem 2.4]{B86}, one can find two weaker assumptions on $B$ that together guarantee existence and uniqueness of the solution of the system \eqref{sXs} and, hence, the similarity of the operators $A - B$ and $ A - X_{11}^*- X_{22}^*$. One of those weaker assumptions appears in Theorem \ref{mt6} below and guarantees the solvability of the first two equations in \eqref{sXs}. The other one is similar and applies to the second pair of equations in \eqref{sXs}.
\erem 

%
%
%
%
%

Next, we would like to obtain better estimates of a perturbed simple eigenvalue and the corresponding eigenvector. 
Consider the function $\Psi: H\to H$ given by
\begeq\label{Psi}
\Psi(z) = (b_1S - B_{22}S -\la B_{12}Sz, e\ra S)z+B_{21}e, \ z\in H.
\eq 
The function $\Psi$ is obtained by applying the second equation in \eqref{sXs} to the normalized eigenvector $e$ of $A$ and denoting $z = X_{21}e$ and $b_1 = \la B_{11}e,e\ra$.
If $\Psi$ has a fixed point $y\in H$, one can verify directly that $e^\prime = e-Sy$ is an eigenvector of $A - B$. The  corresponding eigenvalue is 
$\l^\prime = \l - b_1+b_2$, where   $b_2 = \la B_{12}Sy,e\ra$. Indeed, since $y = Qy$ and $B_{12}Sy =
\la B_{12}Sy, e\ra e$, we have
\[
\bs
(A-B)e^\prime &=(A-B)(e-Sy) = \l e - ASy -Be+BSy \\
& = \l e +y -\l Sy -  b_1e- B_{21}e + B_{21}Sy +B_{22}Sy
\\& = \l e^\prime +(b_1S - B_{22}S -\la B_{12}Sy, e\ra S)y -  b_1e  + B_{12}Sy +B_{22}Sy
\\& = \l e^\prime   -  b_1e^\prime +\la B_{12}Sy, e\ra e^\prime =(\l     -  b_1+ b_2) e^\prime =\l^\prime e^\prime.
\end{split}
\]
Observe that $\|e'\| \ge 1$, that is the eigenvector $e'$ is not normalized unless $y=0$. If $f = \frac{e^\prime}{\|e^\prime\|}$ and $\|e-e^\prime\|\le\varepsilon <1$, we, however, easily get
\[
\|e-f\| = \frac1{\|e^\prime\|}\left\|(\|e^\prime\|-1)e+e-f \right\| \le\frac{2\varepsilon}{1-\varepsilon}.
\] 

We shall use the Browder fixed point theorem to show that the function $\Psi$ in \eqref{Psi} has a fixed point for a sufficiently small $B$. The result is the following theorem, where we employ the notation introduced in this subsection.

\bt\label{mt6}
Suppose that the operator $A$   satisfies Assumption \ref{assa}, and $\l$ is a simple eigenvalue of $A$ with the corresponding eigenvector $e$ of norm $1$. 
Suppose also that the operator $B\in B(H)$ satisfies $\|B_{12}\|\|B_{21}\| \neq 0$ and
\begeq\label{condvec}
\bs
\|b_1S&-B_{22}S\| +2\sqrt{s\|B_{12}S\|\|B_{21}\|} \le 1.
\end{split}
\eq
Let 
\[
r = \frac{1-\|b_1S-B_{22}S\|-\sqrt{(1-\|b_1S-B_{22}S\|)^2-4s\|B_{12}S\|\|B_{21}\|}}{2s\|B_{12}S\|\|B_{21}\|}.
\]
Then $r> 0$ and  the function $\Psi$ in \eqref{Psi} has a fixed point  $y$  in the  closed ball $\mathscr B \subset H$ of radius $r\|B_{21}\|$ centered at $0$. Consequently, 
$e^\prime = e-Sy$ is an eigenvector of $A - B$   corresponding to the eigenvalue  
$\l^\prime = \l - b_1+b_2$, where $b_1 = \la B_{11}e,e\ra$ and $b_2 = \la B_{12}Sy,e\ra$, and we have
\[
\bs
\|e^\prime - e\| &= \|Sy\| \le sr\|B_{21}\| 
\\&\le \frac{s\|B_{21}\|}{1-\|b_1S-B_{22}S\|}
\left(1+ \frac{s\|B_{12}S\|\|B_{21}\|}{(1-\|b_1S-B_{22}S\|)^2}\right), \quad \mbox{and}
\end{split}
\]
\[
\bs
|b_2| &\le r\|B_{12}S\|\|B_{21}\| \\ 
&\le \frac{\|B_{12}S\|\|B_{21}\|}{1-\|b_1S-B_{22}S\|}
\left(1+ \frac{s\|B_{12}S\|\|B_{21}\|}{(1-\|b_1S-B_{22}S\|)^2}\right).
\end{split}
\]
\et

\bpf
To simplify the exposition, we let $m = \|b_1S-B_{22}S\|$ and $n = s\|B_{12}S\|\|B_{21}\|$.
First, consider the quadratic polynomial
\begeq\label{quad}
q(t) = nt^2+(m-1)t +1.
\eq
Clearly, $r$ is one of its roots, which is 
 real and  positive  due to \eqref{condvec}  and $q(0) > 0$.

Second, we show that $\Psi(\mathscr B) \subseteq \mathscr B$. This follows from
\[
\bs
\|\Psi(z)\| &\le \|(b_1S - B_{22}S - \la B_{12}Sz, e\ra S)z\|+\|B_{21}\|
\\&\le (nr^2+mr+1)\|B_{21}\|
= r \|B_{21}\|,
\end{split}
\]
where the last equality holds because $r$ is the root of $q$ in \eqref{quad}.

Third, we derive that $\Psi$ is non-expansive in $\mathscr B$. This follows from
\[
\bs
\|\Psi(w) &-\Psi(z)\| \le \|(b_1S - B_{22}S)(w-z)\|
\\&+\|\la B_{12}Sz, e\ra Sz-\la B_{12}Sz, e\ra Sw+\la B_{12}Sz, e\ra Sw-\la B_{12}Sw, e\ra Sw\|
\\&\le m\|w-z\| + s\|B_{12}S\|\|z\|\|w-z\|+ s\|B_{12}S\|\|w\|\|w-z\|
\\&\le\left(m +2nr\right)\|w-z\|
=\left(1-\sqrt{(1-m)^2-4n}\right)\|w-z\|\le\|w-z\|.
\end{split}
\]
Thus, Browder's fixed point theorem applies and the result is proved. The final estimates for $\|e^\prime - e\|$ and $|b_2|$ are obtained via the second order Tailor approximation of $r$.
\epf

\brem
Observe that the condition $4s\sqrt2\|B\| < 1$ in Theorem \ref{mt5b} implies \eqref{condvec}, so that the assumptions on $B$ in Theorem \ref{mt6} are, indeed, weaker. In particular, \eqref{condvec} is insufficient to reach Conclusion \ref{con}. Observe also that if the inequality in \eqref{condvec} is strict, then one can use the Banach fixed point theorem and obtain the vector $y$ as a limit of simple iterations: $y_{j+1} = \Psi (y_j)$, $y_0 = 0$.
\erem

\brem
If $B_{21} = 0$, we have $\l^\prime = \l-b_1$ and $e^\prime = e$. If $B_{12} = 0$, we also have $\l^\prime = \l-b_1$; for the eigenvector, we then have
$e^\prime = e- S(I-b_1S+B_{22}S)^{-1}B_{21}e$, provided that the inverse is well defined. The latter is true, for example, if \eqref{condvec} holds with a strict inequality. 
\erem

\section{Examples}\label{examples}

In this section, we illustrate our main results with a few examples involving differential and integro-differential operators. 

The Hilbert space $H$ will be the space $L^2 = L^2(\Omega)$ of all (equivalence classes of) Lebesgue square integrable complex-valued functions on a set $\Omega\in\{[0,1], [0,1]^2\}$. To specify the domains of  operators  we will also use the Sobolev spaces $W_2^k = W_2^k(\Omega)$, $k=1,2$, defined in the standard way. In particular, $W_2^1$ is the space of absolutely continuous  functions with the first derivative in $L^2$ and $W_2^2$ is the space of absolutely continuous functions with the first derivative in $W^1_2$. 

\subsection{First order differential operator with an integral perturbation.}\label{newex}

We begin with an example involving integro-differential operators. Such operators have not been analyzed before using the method of similar operators. 

We let $A = \frac d{dt}$ with $D(A) = \{x\in W_2^1([0,1]): x(0) = x(1)\}$. This operator is easily seen to satisfy Assumption \ref{assa}. In particular, we have $\s(A) = \{\l_k = 2\pi ik: k\in\ZZ\}$, yielding $\eta = \frac1{12}$ in \eqref{specAcond} and $\delta  = s^{-1}= 2\pi$ in \eqref{speccond} and \eqref{nG}. The corresponding eigenvectors $e_k = e^{2\pi i k(\cdot)}$ form an orthonormal basis so that the eigenvalues are simple and the resolution of the identity $\mathcal P$ is formed by the rank-one projections $P_kx = \la x, e_k\ra e_k$, $k\in\ZZ$. As always in this paper, the matrix of $A$ is a bi-infinite diagonal operator matrix with non-zero elements given by $\l_kP_k$. 

The operator $B$ is taken to be integral and is defined via a kernel $K\in L^2([0,1]^2)$:
\[
(Bx)(t) = \int_0^1 K(s,t)x(s)ds.
\]
Clearly, $B\in\mathfrak S_2(H)$ and $\|B\|_2 = \|K\|$. As we mentioned in Remark \ref{compB}, it follows that Assumptions \ref{asB} and \ref{asBc} automatically hold. Computing the matrix of $B$, we get that $B_{mn} x = \widehat K(-m, n)\la x, e_n\ra e_m$, $x\in H$, where
\[
\widehat K(m, n) = \int_0^1\int_0^1 K(t,s)e_m(-t)e_n(-s)dtds
\]
 are the Fourier coefficients of the kernel $K$.

Thus, Theorem \ref{mt1} applies if $\|K\| < \frac\pi2$. It yields the similarity of $A - B$ and $A-V$ with $V\in \mathfrak S_2(H)$ having a diagonal matrix.
This gives us an approximation of $\s(A-B)$ up to an $\ell^2$ sequence. In view of Theorem \ref{mt12}, we, in fact, have a weighted $\ell^2$ approximation with the weight $w$ given by $w(k) = (\a_k(B))^{-2}$. More precisely, we have $\sum_k |b_k|^2 w_k< \infty$ for $b_k$ in \eqref{spec1}.

In case Theorem  \ref{mt6} also applies, for example when $\|K\| < \frac\pi{2\sqrt2}$, the approximation can be even further enhanced. We  show it in the following concrete example where we chose a specific kernel for the  perturbation $B$.

\bex
 We let $K(s,t) = s+t$. In this case, we have
$\|B\| \le \|B\|_2 = \|K\| =\left(\int_0^1\int^1_0 (s+t)^2dsdt\right)^{\frac12} =\sqrt{\frac76}$ and
\[
\widehat K(m,n) = 
\left\{\begin{array}{ll}
1,& m=n=0;\\
-\frac1{2\pi i m}, & m\neq 0, n=0;\\
-\frac1{2\pi i n}, & n\neq 0, m=0;\\
0, & mn\neq 0.
\end{array}
\right.
\]

Since $\sqrt{\frac76} <\frac\pi{2}$, Theorem \ref{mt1} applies. Computing 
\[
\a_n (B)\le \left(\frac3{14\pi^2(n-1)}\right)^{\frac14} ,\ n>1,
\]
we get $\sum_k |b_k|^2 \sqrt {|k|}< \infty$ for $b_k$ in \eqref{spec1} via Theorem \ref{mt12}. One, however, gets a much better estimate using the spectral splitting.

Since  $\sqrt{\frac76} <\frac\pi{2\sqrt2}$, Theorems \ref{mt5b} and \ref{mt6} apply for any $\l\in\s(A)$. Let us, first,  use Theorem \ref{mt6} with $\l = 2\pi i k \neq 0$. Observe, that with the notation of the theorem, we have
\[
b_1 = 0,\ \|B_{21}\| = \frac1{2\pi|k|},\   \|B_{12}S\| = \frac1{4\pi^2k^2}, \ \mbox{and}
\]
\[
\frac1{2\pi |k|}\le \|B_{22}S\| \le \|B\|_2 \|S\| = \frac1{2\pi}\sqrt{\frac76}.\footnote{It is clear from these estimates that there is $\kappa > 1$ such that Theorem \ref{mt6} applies for $A-\kappa B$ but Theorem \ref{mt5b} does not.}
\]
The estimates in Theorem \ref{mt6} become 
\[
\|e - e^\prime\| \le\frac{1}{2\pi\left(2\pi|k|-1\right)}\left(1+
\frac{1}{4\pi^2|k|\left(2\pi|k|-1\right)^2}\right)
\]
and
\[
|\l^\prime_k - 2\pi ik| =|b_2| \le\frac{1}{4\pi^2k^2\left(2\pi|k|-1\right)}\left(1+
\frac{1}{4\pi^2|k|\left(2\pi|k|-1\right)^2}\right).
\]
In particular, it follows that we can compute $\s(A-B)$ up to an $O(|k|^{-3})$ sequence.

For the case of $\l = 0$, we have
\[
b_1 = 1,\ \|B_{21}\| = \frac1{2\pi},\  B_{22} = 0, \ \mbox{and } 
\]
\[
\|B_{12}S\| = \left(\sum_{\ell\neq 0} \frac 1{|2\pi i \ell|^4}\right)^{\frac12}= \frac1{4\pi^2} \sqrt{\frac{2\pi^4}{90}}= \frac1{12\sqrt5}
\]
The estimates in Theorem \ref{mt6} become 
\[
\|e - e^\prime\| \le\frac{3\sqrt5}{\pi}\left(2\pi - 1-
\sqrt{(2\pi-1)^2-\frac{1}{3\sqrt5}}\right)\approx 0.0302
\]
and
\[
|\l^\prime_0 - 1| = |b_2| \le\frac{1}{2}\left(2\pi - 1-
\sqrt{(2\pi-1)^2-\frac{1}{3\sqrt5}}\right)\approx0.0071.
\]
\eex

\subsection{Other examples.}

In this subsection, we mention a few examples of operators that have already been analyzed by the method of similar operators. We show that the theory developed here also applies and could be used to streamline the pre-existing arguments. 

\subsubsection{Perturbations with an involution.}

Here, we again have $A = \frac d{dt}$  as in Section \ref{newex}. This time, however, we allow  slightly more flexibility in choosing the domain by letting
\begeq\label{dag}
D(A) = \{x\in W_2^1([0,1]): x(0) = e^{i\pi\theta}x(1)\}
\eq
for some fixed $\theta \in [0,2)$. In the previous example, we restricted ourselves to the choice of $\theta = 0$ to simplify the exposition. We now have that $\s(A)$ consists of simple eigenvalues
$\l_k = \pi i(2k-\theta)$, and the corresponding eigenvectors are $e_k = e^{\l_k(\cdot)}$, $k\in\ZZ$.
Thus, we still have $\eta = \frac1{12}$ in \eqref{specAcond} and $\delta  = s^{-1}= 2\pi$ in \eqref{speccond} and \eqref{nG}. Therefore, the operator $A$   satisfies Assumption \ref{assa}.

The operator $B$ is chosen as in \cite{BKR17, BKU18}. In particular, we let $(Bx)(t) = v(t)x(1-t)$ for some $v\in L^2(0,1) = H$. We cite, for example, \cite{BKh14, SS12} for other recent work where perturbations with an involution were studied.

In this case, we can only guarantee that $B\in \mathfrak L_A(H)$. Thus, we need to make sure that $B$ satisfies Assumptions \ref{asB} and \ref{asBc}. Observe that the matrix elements of $B$ satisfy
\[
\|B_{mn}\|_2 = \left|\int_0^1v(t)e^{\l_n(1-t)}e^{-\l_mt}dt\right| = \left|e^{-\pi i\theta}\widehat{\widetilde v}(m+n)\right| = \left|\widehat{\widetilde v}(m+n)\right|, 
\]
where $\widetilde v = ve^{2\pi i \theta (\cdot)}$ is again in $H$. Thus, Assumption \ref{asBc} clearly holds. Condition \eqref{cGB} in Assumption \ref{asB} is also immediate. Condition \eqref{cBGB} is very difficult to verify directly. It does, however, follow from an integral representation \cite[(2.14)]{BKR17} of the operator $B\Gamma B$ for the case $\theta = 0$. The representation implies immediately that for any choice of $v\in H$ we have $B\Gamma B\in\mathfrak S_2(H)$, when $\theta = 0$. 
In fact, for any $v\in H$, we have
\begeq\label{ineq}
\bs
\frac1{4\pi^2}&\sum_{m,n\in\ZZ}\left|\sum_{\ell\in\ZZ\setminus\{n\}}\!\!\frac{\widehat v(l+m)\widehat v(\ell+n)}{\ell-n}\right|^2 \!= \\&
\int_0^1\!\int_0^1\left|(s-t-\frac12)v(1-s-t)v(2s)\right|^2\!\!dsdt \le\frac94\|v\|^4 <\infty,
\end{split}
\eq
which, in turn, implies
Condition \eqref{cBGB} for any choice of $\theta\in[0,2)$ and $v\in H$. Property \ref{assb3} of Assumption \ref{asB} is also proved the same way as in
\cite{BKR17}.

It follows that Theorem \ref{mt3} and Corollary \ref{cor1d} apply for the operator $A-B$. In particular for the sequences $p = (p_n)$ and $q = (q_n)$ defined in the corollary we have
\[
p_n = e^{-\pi i \theta}\widehat{\widetilde v}(2n) \mbox{ and } 
q_n = \sum_{\ell\in\ZZ\setminus\{n\}}\frac{e^{-2\pi i \theta}(\widehat{\widetilde v}(l+n))^2}{2\pi i(\ell-n)}, \ n\in\ZZ.
\] 
To sum up, we get that the eigenvalues $\l^\prime_n$ of $A-B$ satisfy
\[
\l^\prime_n = \l_n-e^{-\pi i \theta}\widehat{\widetilde v}(2n)-\sum_{\ell\in\ZZ\setminus\{n\}}\frac{e^{-2\pi i \theta}(\widehat{\widetilde v}(l+n))^2}{2\pi i(\ell-n)}-c_n, \ n\in\ZZ,
\]
where $c=(c_n)\in\ell^1$. 

\subsubsection{Dirac operators.}

In this subsection, we show that the analysis of the Dirac operators in \cite{BDS11, BKU19} also falls into the framework developed in this paper. Here, we have $\Omega = [0,1]^2$ and let $A: D(A)\subset L^2(\Omega)\to L^2(\Omega)$ be defined by
\begin{equation}\label{dir1}
(Ay)(t)=i
\begin{pmatrix}
1 & 0 \\
0 & -1
\end{pmatrix}
\frac{dy}{dt},
\end{equation}
where $t\in [0, 1]$ and $D(A) = \{y\in W_2^1([0,1]): y(0) = y(1)\}$. Boundary conditions \eqref{dag} can also be considered but we choose not to burden the exposition.
We then have $\sigma(A)=\left\{{2\pi n}: n\in\mathbb{Z}\right\}$, where each $\lambda_n={2\pi n}$, 
$n\in\mathbb{Z}$, is an eigenvalue of multiplicity two. Moreover, the corresponding eigenspace is given by $\mathrm{span}\,\{e_n^1, e_n^2\}$, where
$$
e_n^1=
\begin{pmatrix}
e_{-n} \\
0
\end{pmatrix},
\quad
e_n^2 =
\begin{pmatrix}
0 \\
e_n
\end{pmatrix},
\quad e_n(t)
=e^{i\lambda_nt}, \quad t\in [0, 1].
$$
Thus, we once again have $\eta = \frac1{12}$ in \eqref{specAcond} and $\delta  = s^{-1}= 2\pi$ in \eqref{speccond} and \eqref{nG}, so that the operator $A$   satisfies Assumption \ref{assa}.

For the perturbation $B$, we let
\begin{equation}\label{dir2}
(By)(t)= V(t)y(t),
\end{equation}
where $t\in [0, 1]$ and 
\begin{equation}\label{dir22}
V(t)=
\begin{pmatrix}
v_1(t) & v_2(t) \\
v_3(t) & v_4(t)
\end{pmatrix},
\end{equation}
$v_j\in L^2[0, 1]$, $1\le  j\le4$.

It is shown in \cite{BKU19} that the operator $A-B$ is similar to $ A - \widetilde B$, where
\begin{equation}\label{dir12}
(\widetilde{B}y)(t)= \begin{pmatrix}
\widehat{v}_1(0) & 0 \\
0 & \widehat{v}_4(0)
\end{pmatrix}
y(t)+
\begin{pmatrix}
0 & u_2(t) \\
u_3(t) & 0
\end{pmatrix}
y(t),
\end{equation}
with $u_2(t) = v_2(t)e^{ig(t)}$, $u_3(t) = v_3(t)e^{-ig(t)}$ and
\[
g(t) = -(\widehat v_1(0)+\widehat v_4(0))t+\int_0^t v_1(\tau)+ v_4(\tau)d\tau, \ t\in[0,1].
\]
We write $\widetilde B = \widetilde B_1 +\widetilde B_2$ where the two operators $\widetilde B_1$ and $\widetilde B_2$ correspond to the two summands on the right hand side of \eqref{dir12}. In \cite{BDS11}, it was proved that $J\widetilde B_2$, $\Gamma\widetilde B_2$, and $\widetilde B_2\Gamma\widetilde B_2$  belong to $\mathfrak S_2(H)$. It was also shown there that Assimption \ref{asB}\eqref{assb3}  holds for $\widetilde B_2$. Consequently, $\widetilde B_2$ satisfies Assumption \ref{asB}.  Moreover, we have $J\widetilde B_1 = \widetilde B_1 \in B(H)$ and $\Gamma\widetilde B_1=0$. It follows that $\widetilde B$ satisfies Assumptions \ref{asB} and \ref{asBb}.

Thus, Theorem \ref{1pre} applies for the operator $A-\widetilde B$, allowing us to conclude that it is similar to the operator $A-J\widetilde B - B_0 = A-\widetilde B_1-\widetilde B_0$, where $B_0, \widetilde B_0\in\mathfrak S_2(H)${\footnote{Recall that $J_m\widetilde B$ is a finite-rank perturbation of $J\widetilde B$.}}. It is not immediately clear if Assumption \ref{abba} holds for $A$ and $\widetilde B$, which prevents us from applying Theorem \ref{mt4}. It is, however, clear that $A-\widetilde B_1$  satisfies Assumption \ref{assa}, so that we can apply Theorem \ref{mt12} for an unperturbed operator $A-\widetilde B_1$ and the perturbation $\widetilde B_0$. It follows, that $A-B$ is similar to $A-\widetilde B_1 - \widetilde J_kX_*$ for some $k\in\ZZ_+$, where the transform $\widetilde J_k$ is defined in the usual way starting from the unperturbed operator $A- \widetilde B_1$, and $X_*\in \M^{\widetilde B_0}$ is the fixed point of the corresponding nonlinear function $\Phi$ defined as in \eqref{bask13}. Therefore, the spectrum $\s(A-B)$ can be estimated up to an $\ell^2$ sequence. A more thorough analysis of the resulting estimates was performed in \cite{BKU19}. There, in most cases, an approximation up to an $\ell^1$ sequence was exhibited. In some degenerate cases, the sequence was proved to be in $\ell^{4/3}$.

\subsubsection{Hill differential operators.}
In this example, we discuss the operators studied in \cite{BP17}. Thus, we let $A = -\frac{d^2}{dt^2}$ be the second order differential operator with the domain
\begeq\label{dag2}
D(A) = \{x\in W_2^2([0,1]): x(0) = e^{i\pi\theta}x(1),\ x^\prime(0) = e^{i\pi\theta}x^\prime(1)\}
\eq
for some $\theta\in(0,1)$. We could also consider $\theta\in\{0,1\}$ but that would require a separate treatment. We have $\s(A) = \{\l_n: n\in\ZZ\}$, where each $\l_n = (\pi(2n-\theta))^2$ is an eigenvalue of multiplicity $1$. The corresponding orthonormal basis of eigenfunctions for $H = L^2([0,1])$ is given by \[e_n(t) = e^{i(\pi(2n-\theta))t},\ n \in \ZZ,\ t \in [0,1].\]
In this case, 
it is also easy to see that $A$ satisfies Assumption \ref{assa}.

For $B$, we take $(Bx)(t) = v(t)x(t)$ with $v\in H$. It is shown in \cite{BP17} that the operator $B$ satisfies Assumptions \ref{asB} and \ref{asBc}. Therefore, Theorem \ref{mt3} and Corollary \ref{cor1d} apply. In this case, using the formulas in \cite{BP17}, we get
\[
p_n = \widehat{ v}(0) \mbox{ and } 
q_n = \frac1{4\pi^2}\sum_{\ell\in\ZZ\setminus\{n\}}\frac{\widehat{ v}(\ell-n)\widehat{ v}(n-\ell)}{ (n-\ell)(n+\ell-\theta)}, \ n\in\ZZ,
\] 
for the sequences $p$ and $q$ in Corollary \ref{cor1d}.

In \cite{BP17}, other types of estimates for $\s(A-B)$ were obtained. We also remark that here $\lim_{m\to\infty} \|\Gamma_m\|_{B(\mathfrak S_2(H))} = 0$, and, consequently, the use of $\M^B$ in the proof of Theorem \ref{mt4} for this case can be avoided.

\section{Acknowledgments}
We thank A.~Caragea, F.~Philipp, and D.~Polyakov for stimulating discussions that helped us improve this manuscript.

\bibliographystyle{siam}
\bibliography{../refs}

\end {document}